%% file: header.tex
\documentclass[12pt]{article}
\usepackage{a4mlnew,arthead}
\usepackage{math302e}
\usepackage{ready}
\usepackage{math40develop}
\usepackage{latexsym,amssymb}
\usepackage{local}

\newcommand{\N}{{\Bbb{N}}}
\newcommand{\Z}{{\Bbb{Z}}}
\newcommand{\C}{{\Bbb{C}}}
\newcommand{\R}{{\Bbb{R}}}

\newcommand{\poly}{$\log$--poly\-homo\-geneous}

\begin{document}
\Subjclass{Primary 58G15; Secondary 58G11}
\Keywords{Pseudodifferential operators, noncommutative residue}
\title{On the noncommutative residue for pseudodifferential operators
with \poly\ symbols}
\author{Matthias Lesch}
\puttitle
\newpage

{
\input introduction.tex

}
\input sec1.tex

\input sec2.tex
\begin{footnotesize}

\input bibliography.tex
\end{footnotesize}




\bigskip
\begin{raggedright}
Matthias Lesch\\
Institut f\"ur Mathematik\\
Humboldt--Universit\"at zu Berlin\\
Unter den Linden 6\\
D--10099 Berlin\\[1em]
email: lesch@mathematik.hu-berlin.de\\

\medskip\noindent
preprints available from\\
http://spectrum.mathematik.hu-berlin.de/$\sim$lesch
\end{raggedright}

\end{document}

%% file: introduction.tex
\section{Introduction and summary of the results}

Let $M$ be a compact smooth Riemannian manifold without boundary. We denote
by $\CL^*(M)$ the algebra of classical (1--step polyhomogeneous)
pseudodifferential operators. $\CL^*(M)$ acts naturally as unbounded
operators on the Hilbert space $\rmL^2(M)$ of square integrable functions.
If $m<-\dim M$ then $\CL^m(M)\subset C_1(\rmL^2(M))$, the space of
trace class operators. 

The $\rmL^2$--trace
does not have a continuation as a trace functional on the whole
algebra $\CL^*(M)$. For assume we had a trace $\tau$ on $\CL^*(M)$ that
extends the $\rmL^2$--trace. 
For $r$ large enough we may choose an elliptic operator
$T\in \CL^1(M)\otimes M_r(\C)$ of nonvanishing Fredholm index. Let $S\in
\CL^{-1}(M)\otimes M_r(\C)$ be a pseudodifferential parametrix. Then
$I-ST, I-TS\in \CL^{-\infty}(M)\otimes M_r(\C)$ and we arrive at
the contradiction
\begin{equation}
   0\not=\ind T=\Tr_{\rmL^2}([T,S])=\tau([T,S])=0.
  \mylabel{G1-I.1}
\end{equation}


However, in his seminal papers \cite{Wod:SALI}, \cite{Wod:NRCIF}
{\sc M. Wodzicki} showed that, up to a constant,
the algebra $\CL^*(M)$ has a unique trace which he called the
noncommutative residue. The noncommutative residue was independently
discovered by {\sc V. Guillemin} \cite{Gui:NPWFADE} as a byproduct
of his so--called "soft" proof of the Weyl asymptotic.

A detailed account of the noncommutative
residue was given by {\sc C. Kassel} \cite{Kas:RNC}.
{\sc B. Fedosov} et al. \cite{\Schretal} generalized the noncommutative
residue to the {\sc Boutet de Monvel} algebra on a manifold with boundary.
Furthermore, {\sc E. Schrohe} \cite{Schro:WNRTOAMCS} studied manifolds
with conical singularities.

There are several ways to define the noncommutative residue.
The global definition which shows the intimate relation 
to $\zeta$--functions and heat kernel expansions
is as follows: given $A\in \CL^a(M)$. Then choose any elliptic operator
$P\in \CL^m(M), m>0,$ whose leading symbol is positive. Then
\begin{eqnarray}
  \Res(A)&:=& m \Res_{s=0} \Tr(A P^{-s})\nonumber\\
   &=&-m \times 
  \;\mbox{coefficient of }\;
      \log t\;\mbox{in the asymptotic}\label{G5-1.2}\\
    &&\mbox{expansion of}\;\Tr(Ae^{-tP})\;\mbox{as}\; t\to 0.\nonumber
\end{eqnarray}
$\Res(A)$ is in fact independent of the $P$ chosen.
This definition uses the fact that the generalized $\zeta$--function
$\Tr(A P^{-s})$ has a meromorphic continuation to $\C$ with simple
poles at $\{\frac{\dim M+a-j}{m}\,|\, j\in \Z_+\}$, $\Z_+:=\{0,1,\ldots\}$.
Via
the Mellin transform
$$\Tr(AP^{-s})=\frac{1}{\Gamma(s)}\int_0^\infty t^{s-1} \Tr(A e^{-tP}) dt,$$
this is (almost) equivalent to the asymptotic expansion
\begin{equation}
\Tr(A e^{-tP})\sim_{t\to 0+}
   \sum_{j=0}^\infty (c_j+c_j'\log t)\, t^{\frac{j-a-\dim M}{m}}
      + \sum_{j=0}^\infty d_j\, t^{j},
   \mylabel{G1-I.2}
\end{equation}
where $c_j'=0$ if $\frac{j-a-\dim M}{m}\not\in \Z_+$
(cf. \cite[Thm 2.7]{GruSee:WPPOAPSBP}).

Another approach to the noncommutative residue was also
discovered independently
by {\sc M. Wodzicki} and {\sc V. Guillemin}.
This approach works
in the framework of symplectic cones. It shows that there is a local formula
\begin{equation}
\Res(A)=(2\pi)^{-\dim M} \int_{S^*M} a(x,\xi)_{-\dim M} |d\xi dx|
\mylabel{G2-I.3}
\end{equation}
for the noncommuative residue in terms of the complete symbol of $A$.
It is a remarkable fact that, although the complete symbol does not
have an invariant meaning, the right hand side of \myref{G2-I.3} is
well--defined.

Asymptotic expansions like \myref{G1-I.2} are an essential feature in the
study of elliptic operators. They can also be achieved for certain operators
with singularities, for example conical singularities \cite{Che:SGSRS},
\cite{BruSee:RESORSO}, \cite{BruLes:SGAC} or boundary value problems
\cite{GruSee:ZEFAPSO}, \cite{GruSee:WPPOAPSBP}, \cite{BruLes:EICNLBVP}. In all these situations no higher
$\log t$ powers occur. However, in \cite{BruHei:AEMPEC} higher $\log t$ powers
show up in equivariant heat trace expansions of the Laplacian.
These higher $\log$--terms are produced by the method of proof and
it is conjectured that in fact all coefficients of
$t^\ga \log^k\,t, k\ge 1,$ vanish.

It was one of our motivations to write this paper to
provide a natural algebra of pseudodifferential operators for which
higher $\log t$--powers occur in the heat expansion \myref{G1-I.2}.
Although our class of operators is not really new, to the best of
our knowledge it was never looked at with regard to noncommutative
residue. Our point of view shows that the existence of the 
noncommutative residue as the unique trace
depends heavily on the absence of higher $\log t$--powers. 

The other motivation arose from the following problem: given $A\in\CL^a(M)$
then by a result of {\sc Wodzicki} $A$ is a sum of commutators
of operators in $\CL^a(M)$ if and only if $\Res(A)=0$. Now, one
can ask whether there is a natural class of pseudodifferential
operators $\subset \rmL^*(M)$ containing the classical ones
such that $A$ is a sum of commutators in this class.

Let us look at the following simple minded analogy between
the noncommutative residue and the ordinary residue for functions.
Consider the algebra $\ca:=\C[z,z^{-1}]$ of {\sc Laurent} polynomials.
Put 
\begin{equation}
    \Res(\sum_{j\in\Z} a_j z^j):=a_{-1}.
\end{equation}
Then the {\sc Laurent} polynomial $f$ has a primitive in $\ca$ if and only
if $\Res(f)=0$. However, if we adjoin the primitive of $z^{-1}$,
$\log z$, we obtain the algebra $\cb:=\ca[\log z]=\C[z,z^{-1},\log z]$
which is filtered by
\begin{equation}
\cb^k:=\{\sum_{j=0}^k f_j \log^j z\,|\, f_j\in \ca\}.
\end{equation}
We define the "higher" residue
\begin{equation}
   \Res_k:\cb^k\longrightarrow \C, \quad \sum_{j=0}^k f_j \log^j z\mapsto
    \Res(f_k).
\end{equation}
Then $f\in\cb^k$ has a primitive in $\cb^k$ if and only if
$\Res_k(f)=0$. However, $f\in \cb^k$ always has a primitive in
$\cb^{k+1}$. This picture can be transfered to the noncommutative
residue for pseudodifferential operators.

The idea is quite simple: instead of looking at classical pseudodifferential
operators we consider pseudodifferential operators whose symbol
have an asymptotic expansion
\begin{equation}
a\sim \sum_{j=0}^\infty a_{m-j},
   \label{G1-I.4}
\end{equation}
where $a_{m-j}$ is {\em \poly }  i.e. 
\begin{equation}
  a_{m-j}(x,\xi)=\sum_{l=0}^k a_{m-j,l}(x,\xi)\log^l|\xi|,
  \label{G1-I.5}
\end{equation}
and $a_{m-j,l}(x,\xi)$ is homogeneous of degree $m-j$ (Definition
\plref{S1-3.1}).
We denote the class of pseudodifferential operators having this
symbol expansion by $\CL^{m,k}(M)$, where $m$ denotes the order
and $k$ the highest $\log$--power occuring in the symbol expansion
(Definition \plref{S1-3.2}).

In fact, this class of operators is not new. It was considered before
by {\sc Schrohe} in his thesis \cite{Schro:CPEPO} where he constructed
the complex powers for elliptic operators in this class.

We had some difficulties to find an appropriate name for
the functions of the form \myref{G1-I.5}. First we tried
"polylogarithmic" but then we were informed by several people
that the polylogarithm has a completely different meaning.
{\sc R. Melrose} suggested to us to use "polyhomogeneous" instead.
But this would conflict with the use of this word in 
\cite{GruSee:WPPOAPSBP}. So we took \poly\ as a compromise.
But still, we also find this term only suboptimal.


In the present paper we show the heat expansion $\Tr(Ae^{-tP})$ for 
$A$ with \poly\  symbol and classical $P$. We use the
method of {\sc Grubb} and {\sc Seeley} \cite{\GS}. In contrast
to \myref{G1-I.2} there occur higher $\log t$--powers
(Theorem \plref{S1-3.5}). 

As a consequence of the heat expansion for $A\in \CL^{m,k}(M)$ we can
define the
"higher" noncommutative residue, $\Res_k(A)$,
 as the coefficient of
the highest $\log t$--power in the expansion of $\Tr(A e^{-tP})$
(Definition \plref{S1-4.1}). It turns
out that this residue has similar properties as the noncommutative
residue of {\sc Wodzicki}, in particular it is independent of the
$P$ chosen. It vanishes on appropriate commutators, i.e.
$$\Res_{k+l}([A,B])=0$$
if $A\in \CL^{m,k}(M)$ and $B\in \CL^{n,l}(M)$
(Theorem \plref{S1-4.4}). There is also a local
formula for $\Res_k$ (Corollary \plref{S1-4.8}). 

In the context of spectral triples higher noncommutative residues
were also discovered and investigated by {\sc A. Connes} and
{\sc H. Moscovici} \cite[Chap. II]{ConMos:LITNG}. In fact,
our algebra of pseudodifferential operators provides examples
of spectral triples with a discrete dimension spectrum
of infinite multiplicity.

As for the {\sc Wodzicki} residue $\Res_k$ is an obstruction
for being a sum of commutators. Namely, if $A\in \CL^{m,k}(M)$ then
there exist $P_1,\ldots, P_N\in \CL^{1,0}(M)$ and $Q_1,\ldots, Q_N
\in \CL^{m,k}(M)$ such that $A-\DST\sum_{j=1}^N [P_j,Q_j]$ is
smoothing if and only if $\Res_k(A)=0$ (Proposition \plref{S1-4.7}).
By a result of {\sc Wodzicki} \cite{Wod:SALI}
(cf. also \cite{Gui:RTCAFIO} for a generalization to Fourier integral
operators) any smoothing operator is in fact
a sum of commutators of classical pseudodifferential operators.
Hence $A\in \CL^{m,k}(M)$ is a sum of commutators if and only if
$\Res_k(A)=0$ (Proposition \plref{S2-4.9}).
However,
$\Res_{k+1}(A)$ is always zero for $A\in \CL^{m,k}(M)$, 
such that
$A$ can always be written as a sum of commutators if one increases
the $\log$--degree by one. As a consequence, there does not exist
any trace functional on the algebra $\CL^{*,*}(M)$ (Corollary \plref{S1-4.9},
but see \myref{G5-5.1}ff.).

For proving the result about commutators we generalize a result due to
{\sc Guillemin} \cite{Gui:NPWFADE} about homogeneous functions on symplectic
cones. Namely, we generalize the notion of symplectic residue, which
is closely related to the noncommuative residue, to functions on
symplectic cones of the form
\begin{equation}
\sum_{j=0}^k f_j \log^j p,
  \mylabel{G1-I.3}
\end{equation}
where $f_j$ is homogeneous and $p$ is positive and homogeneous of degree $0$.
We call these functions \poly. 

{\sc Guillemin}'s result says that a homogeneous function on a symplectic
cone is a sum of Poisson brackets of homogeneous functions if and only
if its symplectic residue vanishes. We define a residue for functions
of the form \myref{G1-I.3} and prove an analogue of {\sc Guillemin}'s
result for these functions.

Finally we generalize the {\sc Kontsevich--Vishik}
trace functional. At the beginning of this section we remarked that
the $\rmL^2$--trace does not have an extension as a trace functional
to classical pseudodifferential operators. The proof uses integer 
order operators. Furthermore, the (higher) noncommutative residues,
which are in some sense the obstructions against the extendability of
the $\rmL^2$--trace as a trace, 
are nontrivial only for integer order operators. 

This gives some evidence that the $\rmL^2$--trace has an extension
to non--integer order operators. Indeed, this is true and the corresponding
functional for classical pseudodifferential operators was discovered
by {\sc Kontsevich} and {\sc Vishik} \cite{KonVis:DEPO,KonVis:GDEO}.
The proof of loc. cit. uses the theory of homogeneous distributions.

We present two alternative approaches to the {\sc 
Kontsevich--Vishik}
trace in the generalized context of
our algebra $\CL^{*,*}$. The first one is completely analogous
to the definition \myref{G5-1.2} of the noncommutative residue
\myref{G5-5.1}. This definition however does not show that
the {\sc Kontsevich--Vishik} trace is given by integration of
a canonical density. This canonical density is constructed in our
second approach. If a pseudodifferential operator is locally given
by
\begin{equation}
  Au(x)=(2\pi)^{-n} \int_{\R^n} \int_U a(x,y,\xi) u(y) e^{i<x-y,\xi>}
      dy d\xi,\label{G5-1.8}
\end{equation}
then it is natural to try (we are bit sloppy with notation here)
\begin{equation}
\Tr(A)=(2\pi)^{-n}\int_M \int_{\R^n} a(x,x,\xi)d\xi dx,\label{G5-1.9}
\end{equation}
and the density we are looking for is
\begin{equation}
(2\pi)^{-n}\int_{\R^n} a(x,x,\xi)d\xi |dx|.
\label{G5-1.10}
\end{equation}
Now the integral \myref{G5-1.10} in general only makes sense if
the order of the operator is $<-\dim M$, i.e. if $A$ is trace class.

It is possible to define a regularized integral for
symbol functions of the form \myref{G1-I.5} (cf. \myref{G5-5.11}).
However, to give \myref{G5-1.10} a coordinate invariant meaning
this regularized integral must satisfy the usual transformation
rule, at least with respect to invertible linear maps.
Proposition \plref{S5-5.2} shows that the transformation rule
is true for non--integer order symbols while for integer order
symbols correction terms show up. This explains why the
{\sc Kontsevich--Vishik} trace exists only for non--integer order
operators.

\bigskip
The paper is organized as follows: in Section 2 we present our generalization
of the symplectic residue. It is divided into two parts. The first part
deals with functions on $\R^n\setminus \{0\}$.
Although this case is very elementary, it does not quite fit into the
framework of symplectic cones and therefore we treat it separately.
In this subsection we follow \cite{\Schretal}. Loc. cit. deals with
the noncommutative
residue for manifolds with boundary. The second part of Section
2 deals with the symplectic residue on compact connected symplectic
cones, and we proceed along the lines of \cite{Gui:NPWFADE}.

In Section 3 we introduce the class $\CL^{*,*}(M)$ of pseudodifferential
operators with \poly\  symbol expansion.

Section 4 contains our main results about the higher noncommutative
residues. Finally, the generalization of the {\sc Kontsevich--Vishik} trace
is presented in Section 5.


This work was supported by Deutsche Forschungsgemeinschaft.

%% file: sec1.tex
\newcommand{\smnull}{\setminus \{0\}}
\renewcommand{\vol}{{\rm vol}}

\def\res{{\rm res}}

\def\myleftmark{Lesch:}
\newcommand{\comment}[1]{\relax}
\renewcommand{\tilde}{\widetilde}

\section{\poly\ functions}
\subsection{\poly\ functions on $\R^n\setminus \{0\}$}

In this section we follow in part \cite[Sec. 1]{\Schretal}.

\begin{dfn}{res-S1.1} A function $f\in\cinf{\R^n\smnull}$ is
called {\em \poly } of degree $(a,k)$ if
\begin{equation}
    f(x)=\sum_{j=0}^l g_j(x) P_j(\log h_j(x)),
  \mylabel{res-G1.1}
\end{equation}
where the $P_j\in \C[t]$ are polynomials, $\deg P_j\le k$, and
$g_j, h_j\in\cinf{\R^n\smnull}$, $h_j>0$, are homogeneous functions,
$g_j$ of degree $a$ and $h_j$ of degree $b_j$.
\end{dfn}

We denote the set of all \poly\ functions of degree
$(a,k)$ by $\cp^{a,k}=\cp^{a,k}(\R^n)$ and put
\begin{equation}
          \cp:= \bigoplus_{a\in\C, k\in\Z_+} \cp^{a,k}.
         \mylabel{res-G1.2}
\end{equation}

More generally, if $M$ is a manifold we denote by $\cp^{a,k}(M,\R^n)$ the
set of functions $f\in \cinf{M\times (\R^n\setminus\{0\})}$ such that 
for each $x\in M$ we have $f(x,\cdot)\in\cp^{a,k}$.

\begin{lemma}{res-S1.2} Each $f\in\cp^{a,k}$ has a unique representation
$$f(x)=\sum_{j=0}^k f_j(x) \log^j |x|$$
with $f_j\in\cp^{a,0}$.
\end{lemma}
In the sequel for $f\in\cp^{a,k}$ the coefficient of $\log^j |x|$
will always be denoted by $f_j$.
\proof The uniqueness is clear. It suffices to prove the existence for
$$f(x)=g(x) \log^l h(x)$$
with $g\in\cp^{a,0}, h\in\cp^{b,0}, l\le k$. Then we write
$$h(x)= |x|^b h(x/|x|)=: |x|^b h_1(x),$$
where $h_1\in\cp^{0,0}$. Thus $\log h_1\in\cp^{0,0}$ and
hence
$$ g(x) \log^l h(x)=\sum_{j=0}^l {\TST {l\choose j}} g(x) (\log^{l-j} h_1(x))
   b^j \log^j |x|.\eqno{\Box}$$

\medbreak

An immediate consequence of Lemma \plref{res-S1.2} is the inclusion
\begin{equation}
    \cp^{a,k}\cp^{b,l}\subset \cp^{a+b,k+l},
    \mylabel{res-G1.5}
\end{equation}
i.e. $\cp$ is a bigraded algebra.

There is an analogue of Euler's theorem for \poly\
functions. Namely, consider
$$f(x) =g(x) \log^l |x|,\quad g\in \cp^{a,0}.$$
Then
\begin{equation}
       \sum_{j=1}^n \frac{\pl f}{\pl x_j} x_j=\frac{d}{d t}\big|_{t=1}
              f(tx)=a f(x) +l g(x) \log^{l-1} |x|,
         \mylabel{res-G1.3}
\end{equation}
from which one derives the identity
\begin{equation}
          \sum_{j=1}^n \frac{\pl}{\pl x_j} (x_j F(x))=f(x),\quad 
  \mbox{\rm if}\quad a\not=-n,
\end{equation}
where
\begin{equation}
         F(x)=g(x) \sum_{j=0}^l \frac{(-1)^{l-j}l!}{j! (n+a)^{l-j+1}}
              \log^j |x|,\quad a\not=-n.
\end{equation}
Thus we have proved
\begin{lemma}{res-S1.3} If $a\not=-n$ then for $f\in\cp^{a,k}$ there
exists a $F\in\cp^{a,k}$ such that
$$ f=\sum_{j=1}^n \frac{\pl}{\pl x_j} (x_j F).$$
\end{lemma}

We turn to the case $a=-n$ which is slightly more subtle.
We denote by
\begin{equation}
       \Delta:=-\sum_{j=1}^n \frac{\pl^2}{\pl x_j^2}
\end{equation}
the (positive) Laplacian in $\R^n$. 

\begin{dfn}{res-S1.4}  Let $f\in\cp^{-n,k}$,
$$  f(x)=\sum_{j=0}^k f_j(x) \log^j|x|.$$
Then we put
$$\res_j(f):=\int_{S^{n-1}} f_j(x) d\vol_S(x),$$
where $\vol_S$ denotes the volume form with respect to the standard
metric on $S^{n-1}$.
\end{dfn}

\begin{lemma}{res-S1.5} Let
\begin{equation}
       f(x)=\sum_{j=0}^k f_j(x) \log^j |x|\in \cp^{-n,k}.
     \mylabel{res-G1.4}
\end{equation}
If $n\not=2$ then there exists a $F\in\cp^{2-n,k}$ with $\Delta F=f$
if and only if $\res_k(f)=0$.

If $n=2$ then there exists a $F\in\cp^{0,k}$ with $\Delta F=f$
if and only if $\res_k(f)=\res_{k-1}(f)=0$.
\end{lemma}
\proof We identify $\R^n\smnull$ with $\R_+\times S^{n-1}$ via
$$\R_+\times S^{n-1}\to \R^n\smnull, \quad (r,\theta)\mapsto r\theta.$$
Then we have
$$ f(r\theta)=\sum_{j=0}^k g_j(\theta) r^{-n} \log^j r.$$
The Laplacian is given by
$$\Delta= -r^{1-n}\frac{\pl}{\pl r} (r^{n-1}\frac{\pl}{\pl r})+
           r^{-2}\Delta_S.$$
For $F$ we make an Ansatz
$$F(r\theta)= \sum_{j=0}^k F_j(\theta) r^{2-n} \log^j r$$
and applying $\Delta$ to $F$ we obtain the system of equations
\alpheqn
\begin{eqnarray}
         \Delta_S F_k&=&g_k,\mylabel{res-G1.8a}  \\
         \Delta_S F_{k-1}&=& g_{k-1} -k(n-2) F_k,\mylabel{res-G1.8b}\\
         \Delta_S F_j&=& g_j-(j+1)(n-2) F_{j+1} +(j+1)(j+2)F_{j+2},\quad
               j\le k-2.\mylabel{res-G1.8c}
\end{eqnarray}
\reseteqn
The equation $\Delta_S F_k=g_k$ has a solution if and only if
$g_k$ is orthogonal to $\ker \Delta_S$ which consists of the constants.
Hence the first equation has a solution if and only if $\res_k(f)=0$.
If $n=2$ by the same reasoning the second equation has a solution
if and only if $\res_{k-1}(f)=0$. Hence we have proved the 'only if'
part of the assertion.

Now let $n\not=2$ and assume $\res_k(f)=0$. Let $F_k^0$ be the
unique solution of \myref{res-G1.8a} with $\int_{S^{n-1}} F_k^0=0$. Put
$$F_k:= F_k^0 + \frac{1}{k(n-2)} \int_{S^{n-1}} g_{k-1}(x) d\vol_S(x).$$
Then $\Delta F_k=g_k$ and
$$ k(n-2) \int_{S^{n-1}} F_k(x) d\vol_S(x)=\int_{S^{n-1}} g_{k-1}(x) d\vol_S(x).$$
Hence there exists a unique $F_{k-1}^0$ with $\Delta F_{k-1}^0
=g_{k-1}-k(n-2)F_k$. Proceeding in this way we obtain a solution
of \myref{res-G1.8a}--\myref{res-G1.8c}.

The case $n=2$ is treated similar.
\endproof

\begin{lemma}{res-S1.7} For $f\in\cp^{1-n,k}$ we have
$\res_k(\frac{\pl}{\pl x_j} f)=0$.
\end{lemma}
See \cite[Lemma 1.2]{FedGolLeiSch:NRMB} for another proof of this lemma
in the case $k=0$.
\proof
Let
$$f(x) =\sum_{i=0}^k f_i(x) \log^i|x|,\quad f_i\in\cp^{1-n,0}.$$
Then
$$ \frac{\pl f}{\pl x_j}(x)=\sum_{i=0}^k \frac{\pl f_i}{\pl x_j}(x) \log^i|x|
      +\sum_{i=0}^{k-1} \frac{x_j f_{i+1}(x)}{|x|^2}(i+1)\log^i|x|,$$
hence
$$\res_k(\frac{\pl}{\pl x_j} f)=
    \int_{S^{n-1}} \frac{\pl f_k}{\pl x_j} d\vol_S.$$
We consider the vector field
$$X_j(x):=(1-x_j^2)\frac{\pl}{\pl x_j}-\sum_{i\not=j} x_j x_i \frac{\pl}{\pl x_i}.$$
Obviously,
$X_j|S^{n-1}$ is tangential to $S^{n-1}$. One checks by direct
calculation 
\begin{equation}
  {\rm div}_S X_j=(1-n) x_j.
  \mylabel{res-G1.9}
\end{equation}
On the other hand, by Euler's identity,
\begin{eqnarray*}
X_j f_k&=&\frac{\pl f_k}{\pl x_j}-x_j \sum_{i=1}^n x_i \frac{\pl f_k}{\pl x_i}\\
    &=& \frac{\pl f_k}{\pl x_j}+(n-1) x_j f_k,
\end{eqnarray*}
hence in view of \myref{res-G1.9}
\begin{epeqnarray}
\int_{S^{n-1}} \frac{\pl f_k}{\pl x_j}d\vol_S&=&
         \int_{S^{n-1}} X_j f_k +(1-n) x_j f_k d\vol_S\\
      &=& \int_{S^{n-1}} (-{\rm div}_S (X_j)+ (1-n) x_j) f_k d\vol_S=0.
\end{epeqnarray}

\begin{prop}{res-1.6} Let $f\in\cp^{a,k}$. Then there exist
$f_j\in \cp^{a,k}$ such that $f=\sum_{j=1}^n \frac{\pl}{\pl x_j} f_j$
if and only if $a\not=-n$ or $\res_k(f)=0$.
\end{prop}
\proof The case $a\not=-n$ follows from Lemma \plref{res-S1.3}.
Let $a=-n$. If $\res_k(f)=\res_{k-1}(f)=0$ then by Lemma \plref{res-S1.4}
there exists $F\in\cp^{2-n,k}$ such that $f=\Delta F$ and we reach
the conclusion in this case.

Thus it suffices to show that the function
$$  |x|^{-n} \log^{k-1} |x|$$
is a sum of derivatives. But
$$  \sum_{j=1}^n \frac{\pl}{\pl x_j}\big(x_j  |x|^{-n}\log^k|x|)
     = k |x|^{-n} \log^{k-1}|x|$$
does the job.
\endproof

We will also need the results of 
\cite[(1.1)--(1.5), Lemma 1.1 and 1.3]{\Schretal} in our context.
We briefly summarize the facts. Let
\begin{equation}
   \sigma:=\sum_{j=0}^n(-1)^{j+1} \xi_j d\xi_1\wedge \ldots\wedge
    \widehat{d\xi_j}\wedge\ldots\wedge d\xi_n\in \Omega^{n-1}(\R^n).
\end{equation}
$\sigma|S^{n-1}$ is the volume form. Moreover, for $f\in \cp^{-n,0}$ the form
$f\sigma$ is closed by Euler's theorem \myref{res-G1.3}.
Thus for any bounded domain $D\subset \R^n$,
$0\in D$, with smooth boundary we have for $f\in\cp^{-n,0}$
\begin{equation}
   \res_k(f)=\int_{S^{n-1}} f_k(x) d\vol_{S^{n-1}}(x)=\int_{S^{n-1}} f_k\sigma
   =\int_{\pl D} f_k\sigma.
   \mylabel{res-G1.10}						   
\end{equation}
Moreover, if $T\in {\rm GL}(n,\R)$ then
\begin{eqnarray}
   \res_k(T^*f)&=&\int_{S^{n-1}} f_k\circ T \sigma =\sgn(\det T) \int_{T(S^{n-1})}
       f_k T^{-1\, *}\sigma\nonumber\\
    &=&\frac{1}{|\det T|} \int_{T(S^{n-1})} f_k \sigma
      =\frac{1}{|\det T|}\res_k(f).
   \mylabel{res-G1.11}
\end{eqnarray}

Finally we note

\begin{lemma}{res-S1.8}{\rm \cite[Lemma 1.3]{\Schretal}}
Let $f\in \cp^{a,k}$. Then
$\res_k(\xi^\ga \pl^\gb f)=0$ if $|\gb|>|\ga|$.
\end{lemma}
\proof This follows from Lemma \plref{res-S1.7} by induction since if
$\pl^\gb=\frac{\pl}{\pl \xi_j}\pl^\gamma$ then
$$
\xi^\ga \pl^\gb f= \frac{\pl}{\pl \xi_j}(\xi^\ga \pl^\gamma f)
  -\frac{\pl \xi^\ga}{\pl \xi_j} \pl^\gamma f.\epformel
$$

\subsection{\poly\ functions on symplectic cones}

In this section we study \poly\ functions on an arbitrary
symplectic cone. Our exposition parallels \cite[Sec. 6]{Gui:NPWFADE}.

Let $Y$ be a symplectic cone. This is a principal bundle
$$\pi:Y\to X$$
with structure group $\R_+$. Denote by $\varrho_a:Y\to Y$ the action of
$a\in\R_+$. That $Y$ is a symplectic cone means that $Y$ is symplectic,
with symplectic form $\omega$, and
$$\varrho_a^*\go=a \go.$$

We assume furthermore that $Y$ is {\em connected} and $X$ is 
{\em compact}.

The main example of course is the cotangent bundle with the zero
section removed, $T^*M\setminus 0$, over a compact connected manifold
$M$ of dimension $\dim M>1$.

\begin{dfn}{res-S2.1} A function $f\in\cinf{Y}$ is called
{\em \poly } of degree $(a,k)$ if
$$ f=\sum_{j=0}^l g_j P_j(\log h_j)$$
with $g_j, h_j\in\cinf{Y}, P_j\in\C[t]$, where $g_j$ is homogeneous
of degree $a$, $h_j$ is homogeneous of
degree $b_j$, and $h_j>0$ everywhere. Furthermore,
$\deg P_j\le k$.
\end{dfn}

Again, we denote the set of all \poly\ functions of degree
$(a,k)$ by $\cp^{a,k}$. Then \myref{res-G1.2}, \myref{res-G1.5}
hold similarly.

We fix, for once and for all, $p\in\cp^{1,0}$ such that
$p$ is everywhere positive. Then by Euler's identity we have
$dp\not=0$ everywhere. We put
$$Z:=\{y\in Y\,|\, p(y)=1\}.$$
$p$ plays the role of $|\cdot|$ and $Z$ plays the role of $S^{n-1}$
in the preceding section.

\begin{lemma}{res-S2.2} Each $f\in\cp^{a,k}$ has a representation
$$f=\sum_{j=0}^k f_j \log^j p.$$
Furthermore, $f_k$ is independent of the choice of $p$.
\end{lemma}
\proof Consider
$$g\,\log^l h$$
with $g\in\cp^{a,0}, h\in\cp^{b,0}$. Then $h_1:=p^{-b} h$ is
$0$--homogeneous and positive, hence $\log h_1\in\cp^{0,0}$.
Thus
$$ g \,\log^l h=\sum_{j=0}^l {\TST {l\choose j}} g(x) (\log^{l-j} h_1)
   b^j \log^j p,$$
from which we see that the coefficient of $\log^l p$ is independent
of $p$.
\endproof

\begin{dfn}{res-S2.3} We put
$$\res_k(f):=\Res\, f_k,$$
where $\Res\, f_k$ is the symplectic residue of {\sc Guillemin}.
By the preceding Lemma $\res_k$ is well defined.
\end{dfn}

For the convenience of the reader and since we have to
introduce some notation anyway we briefly
recall the definition of the symplectic residue
(cf. \cite[Sec. 6]{Gui:NPWFADE}).

Via $\Phi_t:=\varrho_{e^t}$ we obtain a one parameter group of
diffeomorphisms of $Y$. Let $\Xi\in\cinf{TY}$ be the infinitesimal
generator of this group. Put
$\ga:=i_\Xi\omega\in \Omega^1(Y),$
and let $$\mu:= \ga\wedge \go^{n-1}.$$
If $f\in\cp^{-n,0}$ then the form $f\mu$ is horizontal and invariant,
hence there is a unique $(2n-1)$--form $\mu_f$ such that
$f\mu=\pi^*\mu_f$. Then
$$\Res(f):=\int_X \mu_f.$$
One can show (cf. \cite[Proof of Lemma 6.3]{Gui:NPWFADE}) that also
$$\Res(f)=\int_Z f\mu.$$

We denote by $\{\cdot,\cdot\}$ the Poisson bracket associated with
the symplectic structure.

\begin{lemma}{res-S2.4} If $f\in\cp^{a,k}, g\in\cp^{1,l}$ then
$\{f,g\}\in\cp^{a,k+l}$ and  $\res_{k+l}\{f,g\}=0$.
\end{lemma}
\proof W.l.o.g. we may assume
$$ f=\phi \,\log^k p,\quad g=\psi\, \log^l p,$$
with $\phi\in\cp^{a,0},\psi\in\cp^{1,0}$. Then
\begin{eqnarray*}
      \{f,g\}&=&  \{\phi,\psi\} \log^{k+l}p +k \{p,\psi\} \phi p^{-1}\log^{k+l-1} p\\
           && +\{\phi,p\} l \psi p^{-1}\log^{k+l-1} p,
\end{eqnarray*}
hence $\{f,g\}\in\cp^{a,k+l}$ and $\res_{k+l}\{f,g\}=\Res \{\phi,\psi\}=0$
by \cite[Prop. 6.1]{Gui:NPWFADE}.
\endproof

Now we can state the generalization of \cite[Thm. 6.2]{Gui:NPWFADE}
to \poly\ functions.

\begin{theorem}{res-S2.5} {\rm 1.} If $a\not=-n$ then $\{\cp^{1,0},\cp^{a,k}\}=
\cp^{a,k}. $

{\rm 2.} If $a=-n$ then $\{\cp^{1,0},\cp^{a,k}\}=\ker \res_k\subset
\cp^{a,k}.$
\end{theorem}
\proof We follow the proof of loc. cit. and choose functions
$g_1,\ldots, g_N\in \cp^{1,0}$ such that their differentials span
the cotangent space of $Y$ at every point.
Let
$$D_i :\cp^{a,k}\to \cp^{a,k} , \quad D_i f:= \{g_i,f\}.$$
We introduce a pre--Hilbert space structure on $\cp^{a,k}$ as follows:
we identify $\cp^{a,k}$ with $\cinf{Z,\C^{k+1}}$ via
\begin{equation}
 \sum_{j=0}^k f_j \log^j p\mapsto (f_0,\ldots,f_k)|Z.
  \mylabel{res-G2.4}
\end{equation}
Next let $\nu$ be the restriction of $\mu$ to $Z$. This is a volume form
and hence it defines a $L^2$--structure on $\cinf{Z,\C^{k+1}}\simeq \cp^{a,k}$.

Now consider $f\,\log^j p, f\in \cp^{a,0}$. Then
$$ D_i (f\log^j p)= (D_if) \log^jp+j f p^{-1} (D_i p) \log^{j-1} p.$$
Thus putting
$$q_i:= p^{-1} D_i p\in\cp^{0,0}$$
and $d_i:= D_i|\cp^{a,0}$ the identification \myref{res-G2.4}
transforms $D_i$ into
\begin{equation}
         D_i\cong\left(\matrix{
       d_i    & q_i    & 0      &  \ldots  & 0\cr
        0     & d_i    & 2 q_i  &  \ldots  & 0 \cr
       \vdots & \ddots & \ddots &\ddots    & \vdots\cr
              &        &        & d_i      & k q_i \cr
        0     &  \ldots&        & 0        & d_i\cr}\right)\;:
      \cinf{Z,\C^{k+1}}\to\cinf{Z,\C^{k+1}}.
      \mylabel{res-G2.5}
\end{equation}
Now we consider
\begin{equation}
      \Delta:=\sum_{j=1}^N D_i D_i^t.
\end{equation}
$\Delta$ is a self--adjoint elliptic differential operator on
$\cinf{Z,\C^{k+1}}$. Hence
$$\Im \Delta^\perp=\ker \Delta.$$
Therefore, consider
$$f=\sum_{j=0}^k f_j \log^j p\in\ker \Delta.$$
In view of \myref{res-G2.5}
$$ d_i^t f_j+j q_i f_{j-1}=0, \quad i=1,\ldots, N; \;j=0,\ldots,k,$$
where we have put $f_{-1}:=0$.

By \cite[(6.15)]{Gui:NPWFADE} we have
$$d_i^t=-d_i + (2a+n) q_i.$$
Abbreviating $r:=-(n+2a)$ we obtain
$$-\{g_i,f_j\}-r p^{-1} f_j \{ g_i,p\} +j f_{j-1} p^{-1} \{g_i,p\}=0,$$
or
\begin{equation}
\{g_i,p^r f_j\}=j f_{j-1} p^{r-1} \{g_i,p\}.
   \mylabel{res-G2.6}
\end{equation}
Since the differentials of the $g_i$ span the cotangent space at every
point we conclude that $p^r f_0=c_0$ is constant.

By induction we assume that $p^r f_{j'}=0$ for $j'<j-1$ and
$p^r f_{j-1}=c_{j-1}$ is constant. Then, in view of \myref{res-G2.6}
$$ \{g_i,p^r f_j\}=j c_{j-1} \{g_i,\log p\},$$
thus
$$\{g_i,p^r f_j-j c_{j-1} \log p\}=0,$$
(cf. \cite[Lemma 6.6]{Gui:NPWFADE}).
Again since the differentials of the $g_i$ span the cotangent space this
implies
$$  f_j= c_j p^{-r} +j c_{j-1} p^{-r} \log p.$$
But since $f_j\in\cp^{a,0}$ the constant $c_{j-1}$ must be $0$.

Summing up we have proved
$$f_j=0, j<k,\quad f_k= c p^{-r}.$$
This implies $c=0$ or $r=-a$. Since $r=-(n+2a)$ this is equivalent
to $c=0$ or $r=n$. Since
$\int_Z p^{-r} \nu\not=0$ we reach the conclusion.
\endproof

%% file: sec2.tex
\section{Pseudodifferential operators with \poly\ symbols}

Let $M$ be a smooth manifold of dimension $n$. We denote by
$\rmL^*(M)$ the algebra of pseudodifferential operators with complete
symbols of {\sc H\"ormander} type (1,0). I. e. if $U$ is an open chart
then $A\in \rmL^m(U)$ can be written
\begin{equation}
   Au(x)=(2\pi)^{-n} \int_{\R^n} \int_U a(x,y,\xi) u(y) e^{i<x-y,\xi>}
      dy d\xi,
   \label{G1-3.1}
\end{equation}
where $a\in\cinf{U\times U\times \R^n}$ satisfies
\begin{equation}
   |\pl_x^\ga\pl_y^\gb \pl_\xi^\gamma a(x,y,\xi)|\le 
    C_{\ga,\gb,\gamma,K} (1+|\xi|)^{m-|\gamma|},
   \label{G1-3.2}	 
\end{equation}
for $\xi\in\R^n, (x,y)\in K$, and $K\subset U$ compact.
The class of these symbols is as usual denoted by $\rmS^m(U\times U,\R^n)$.

We denote by $\rmL^*(M,E)$ the algebra of pseudodifferential operators 
acting on sections of the $C^\infty$ vector bundle $E$.

We are now going to introduce a subclass of $\rmL^*(M)$ which generalizes
the classical pseudodifferential operators.

We proceed locally and introduce \poly\ symbols:

\begin{dfn}{S1-3.1} For an open set $U\subset \R^n$ we denote
by $\CS^{m,k}(U,\R^n)$ the set of symbols 
$a\in\cap_{\eps>0} \rmS^{m+\eps}(U,\R^n)$
having an asymptotic expansion
$$
   a\sim \sum_{j=0}^\infty \psi(\xi) a_{m-j}(x,\xi)
$$
where $a_{m-j}\in \cp^{m-j,k}(U,\R^n)$ and $\psi\in\cinf{\R^n}$
with $\psi(\xi)=0$ for $|\xi|\le 1/4$ and $\psi(\xi)=1$ for $|\xi|\ge 1/2$.
\end{dfn}

\begin{remark}
1. $\CS^{m,0}(U,\R^n)$ are just the classical (1-step polyhomogeneous)
symbols.

2.  One could more generally consider symbols having expansions
\begin{equation}
   a\sim \sum_{j=0}^\infty a_j,
   \label{G1-3.3}  
\end{equation} where $a_j\in \cp^{m_j,k}(U,\R^n),
\lim_{j\to\infty} m_j=-\infty$. These symbols were already considered
by {\sc Schrohe} \cite{Schro:CPEPO} who constructed the complex powers for the
corresponding class of elliptic pseudodifferential operators.

Our class of operators, which is more closely to classical operators,
is large enough for our purposes. However, our results (with suitable
modifications) remain true for the slightly more general symbols
\myref{G1-3.3}.
\end{remark}

\begin{dfn}{S1-3.2} We denote by $\CL^{m,k}(U)$ the class of 
pseudodifferential operators which can be written in the form
\myref{G1-3.1} with $a\in \CS^{m,k}(U\times U,\R^n)$.
\end{dfn}


It is fairly straightforward to check that this class of pseudodifferential
operators satisfies the usual rules of calculus. We summarize
the results for the convenience of the reader.

\begin{prop}{S1-3.3} 
{\rm 1.} If $A\in \CL^{m,k}(U)$ is properly supported
then the complete symbol $\sigma_A$ is in $\CS^{m,k}(U,\R^n)$
and
\begin{eqnarray}
     \sigma_A&\sim& \sum_{\ga\in \Z_+^n} \frac{i^{-|\ga|}}{\ga!} \pl_\xi^\ga
      \pl_y^\ga a(x,y,\xi)|_{y=x}\nonumber\\
       &\sim& \sum_{j=0}^\infty \sum_{|\ga|+l=j}  \frac{i^{-|\ga|}}{\ga!}
     \pl_\xi^\ga\pl_y^\ga a_{m-l}(x,y,\xi)|_{y=x}.\mylabel{G1-3.4}
\end{eqnarray}

{\rm 2.} If $A\in \CL^{m,k}(U), B\in \CL^{m',k'}(U)$ are properly supported
then\\ $AB\in \CL^{m+m',k+k'}(U)$ and
 \begin{eqnarray}
    \sigma_{AB}&\sim& \sum_{\ga\in \Z_+^n} \frac{i^{-|\ga|}}{\ga!} (\pl_\xi^\ga
     \sigma_A(x,\xi))  \pl_x^\ga \sigma_B(x,\xi)\nonumber\\
      &\sim& \sum_{j=0}^\infty \sum_{|\ga|+l+l'=j}  \frac{i^{-|\ga|}}{\ga!}
    (\pl_\xi^\ga  a_{m-l}(x,\xi))  \pl_x^\ga b_{m'-l'}(x,\xi),\mylabel{G1-3.5}
\end{eqnarray}
if $\sigma_A\sim \sum a_{m-j}, \sigma_B\sim \sum b_{m'-j}$.

{\rm 3. } If $A\in \CL^{m,k}(U)$ then $A^t$ in $\CL^{\overline{m},k}(U)$ and

\begin{eqnarray}
     \sigma_{A^t}&\sim& \sum_{\ga\in \Z_+^n} \frac{i^{-|\ga|}}{\ga!} \pl_\xi^\ga
      \pl_x^\ga \sigma_A(x,\xi)^*\nonumber\\
       &\sim& \sum_{j=0}^\infty \sum_{|\ga|+l=j}  \frac{i^{-|\ga|}}{\ga!}
     \pl_\xi^\ga\pl_x^\ga a_{m-l}(x,\xi)^*.\mylabel{G1-3.6}
\end{eqnarray}
\end{prop}
\proof This follows from \cite[Theorems 3.1--3.4]{Shu:POST} and the
obvious inclusions
$$\pl_\xi^\ga \pl_x^\gb \CS^{m,k}(U)\subset \CS^{m-|\ga|,k}(U),
\CS^{m,k}(U)\cdot \CS^{m',k'}(U)\subset \CS^{m+m',k+k'}(U).\epformel$$

\begin{prop}{S1-3.4} Let $\gk:U\to V$ be a diffeomorphism and
$A\in \CL^{m,k}(U)$ properly supported. Then $\gk_* A$ given
by $(\gk_*A)u:= (A(u\circ \gk))\circ \gk^{-1}$ is in
$\CL^{m,k}(V)$ and
\begin{equation}
    \sigma_{\gk_*A}(y,\eta)|_{y=\gk(x)}\sim
     \sum_{\ga\in \Z_+^n} \frac{1}{\ga!} (\pl_\xi^\ga \sigma_A)(x,
      D\gk^t(x) \eta) D_z^\ga e^{i <\gk_x''(z),\eta>}|_{z=x},
    \mylabel{G1-3.7}
\end{equation}
where 
$$\gk_x''(z)= \gk(z)-\gk(x)- D\gk(x)(z-x).$$

Furthermore,
$$\Phi_\ga(x,\eta):=D_z^\ga e^{i <\gk_x''(z),\eta>}|_{z=x}$$
is a polynomial in $\eta$ of degree not higher than $|\ga|/2$.
\end{prop}
\proof Except the assertion $\gk_*A\in \CL^{m,k}(V)$ this is
\cite[Thm. 4.2]{Shu:POST}. $\gk_*A\in \CL^{m,k}(V)$ now follows
from \myref{G1-3.7}.\endproof

As usual, this result allows to define operators of the class $\CL^{m,k}$ on
manifolds. For a smooth manifold $M$ we denote by $\CL^{m,k}(M)$ the
corresponding space of pseudodifferential operators. We note
that \myref{G1-3.7} shows that the leading symbol of $A\in \CL^{m,k}(M)$
can be considered as an element of $\cp^{m,k}(T^*M)$. Because of its
importance we single out this observation:

\begin{prop}{S2-3.5} The leading symbol which is locally defined by
$\sigma_\rmL^m(A):=a_m$ induces a surjective linear map
$$\sigma_\rmL:\CL^{m,k}(M)\longrightarrow \cp^{m,k}(T^*M)\cong
\cinf{\rmS^*M,\C^{k+1}}$$
with $\ker \sigma_\rmL=\CL^{m-1,k}(M)$. Furthermore, for $A\in \CL^{a,k}(M),
B\in \CL^{b,l}(M)$ we have $\sigma_\rmL^{a+b}(AB)=\sigma_\rmL^a(A)\sigma_\rmL^b(B).$
\end{prop}
\proof That $\sigma_\rmL$ is well defined 
follows immediately from Proposition \plref{S1-3.4}.
Similar to \ref{res-G2.4} the isomorphism 
$\cp^{m,k}(T^*M)\cong\cinf{\rmS^*M,\C^{k+1}}$
is given by
$$\sum_{j=0}^k f_j \log^j |.|\mapsto (f_0,\ldots,f_k)|\rmS^*M,$$
where $\rmS^*M$ denotes the cosphere bundle of $M$.
Surjectivity of $\sigma_\rmL$ is proved by the standard construction
of gluing together local pseudodifferential data.
\endproof

For defining the \poly\ residue of an operator $A\in \CL^{a,k}(M)$
we need the meromorphic continuation of the function $\Tr(A P^{-s})$
for a classical elliptic pseudo\-differential operator $P$. In principle, this
could be done for any elliptic $P\in \CL^{m,k}(M)$. The meromorphic
continuation of $\Tr(P^{-s}), P\in \CL^{m,k}(M)$, was already proved
in \cite{Schro:CPEPO}. However, for treating $\Tr(A P^{-s})$ one
has to modify Sections 1 and 2 of \cite{\GS}. 
Therefore, we decided to content ourselves to classical elliptic
$P\in \CL^{m,0}(M)$, which is enough for our purposes.
Then the method of loc. cit. directly applies. Nevertheless,
as mentioned on p. 488 of loc. cit. the results there could be 
generalized to operators of class $\CL^{m,k}$.

First we state the expansion result for the resolvent (cf.
\cite[Thm. 2.7]{\GS}). For the definition of
ellipticity with parameter see loc. cit. Def. 2.7.

\begin{theorem}{S1-3.5} Let $M$ be a compact manifold of dimension $n$,
$E$ a $C^\infty$ vector bundle over $M$, and $A\in \CL^{a,k}(M,E)$.
Assume furthermore that $P\in \CL^{m,0}(M,E)$ is elliptic with
parameter $\mu\in\Gamma$. Then for $\gl\in -\Gamma^m:=\{ -\mu^m\,|\,
\mu \in \Gamma\}$ and $N\in \N$ with $-Nm+a<-n$ the kernel
$(A(P-\gl)^{-N})(x,y)$ of $A(P-\gl)^{-N}$ satisfies on the diagonal
 $$    (A(P-\gl)^{-N})(x,x)\sim \sum_{j=0}^\infty \sum_{l=0}^{k+1}
       c_{jl}(x) \gl^{\frac{n+a-j}{m}-N} \log^l\gl
     \quad + \sum_{j=0}^\infty d_j(x) \gl^{-j-N},$$
as $|\gl|\to\infty$, uniformly in closed subsectors of $\Gamma$.

Furthermore, $c_{j,k+1}=0$ if $\frac{j-a-n}{m}\not\in\Z_+$.
\end{theorem}
\proof This is proved similar to \cite[Thm 2.7]{\GS}. For the convenience
of the reader we indicate the steps of the proof
of loc. cit. that have to be modified. During this proof we will use
freely the notation of \cite{\GS}.

Let $\cp$ be a pseudodifferential parametrix of $(P-\gl)^{-N}$ as
constructed in \cite[p. 503]{\GS}. Then
\begin{equation}
A(P-\gl)^{-N}=A\cp^N+\sum_{l\ge 1} A\cp R_{(l)}
\end{equation}
where $\cp R_{(l)}\in{\rm OP}(S^{-\infty,-(l+1)m})$ and thus
also $A\cp R_{(l)}\in {\rm OP}(S^{-\infty,-(l+1)m})$. Since the symbol
expansion of $A\cp R_{(l)}$ depends on $\mu$ only as a rational
function of $\gl=-\mu^m$ the kernel of $A\cp R_{(l)}$ has the
asymptotic expansion
\begin{equation}
K_{A\cp R_{(l)}}(x,x,\gl)\sim_{\gl\to\infty}
    \sum_{\sigma=0}^\infty c_{l,\sigma}(x) \gl^{-l-N-\sigma}
\end{equation}
(\cite[p. 504]{\GS}).

To expand the kernel of $A\cp^N$ we note that
$A\cp^N={\rm OP}(q)$ where $q\in S^{a+\eps-Nm,0}\cap S^{a+\eps,-Nm}$
for every $\eps>0$. Note that although $A$ is not weakly polyhomogeneous
in the sense of \cite{\GS} we have $A\in \cap_{\eps>0} L^{a+\eps}(M)$.

Now $q$ has an expansion $q\sim \sum_{j\ge 0} q_j$, where
\begin{equation}
q_j =\sum_{|\ga|+l+l'=j}\frac{i^{-|\ga|}}{\ga!}
   (\pl_\xi^\ga a_{a-l}) \pl_x^\ga p_{-mN-l'},
\end{equation}
$a_{a-l}\in \cp^{a-l,k}(T^*M)$, $p_{-mN-l'}$ is $-mN-l'$--homogeneous
in $(\xi,\mu)$ for $|\xi|\ge 1$.

The kernels of the remainders $r_J=q-\sum_{0\le j<J} q_j$ are expanded
as in \cite[(2.3)]{\GS} which gives again integer powers $\gl^{-N-l}$
in the expansion of the kernel. Here we have to note again that
$p_{-mN-l'}(x,\xi,\mu)$ is a rational function in $-\mu^m$.

Picking one of the summands of $q_j$ we are finally facing the
problem of expanding the integral
\begin{equation}
\int_{\R^n} a(x,\xi)\,\log^t|\xi| \,b(x,\xi,\mu)\,d\xi,\quad t\le k,
  \mylabel{G3-4.1}
\end{equation}
where $a(x,\xi)$ is $(a-l-|\ga|)$--homogeneous for $|\xi|\ge 1$ and 
$b(x,\xi,\mu)$ is $(-mN-l')$--homogeneous in $(\xi,\mu)$ for $|\xi|\ge 1$.
Furthermore, $b(x,\xi,\mu)$ is a rational function of $-\mu^m$ with
coefficients homogeneous in $\xi$.

Now we split the integral \myref{G3-4.1} into the three pieces
$|\xi|\ge|\mu|, 1\le |\xi|\le |\mu|, |\xi|\le 1$ as in 
\cite[(2.6)]{\GS}.

By homogeneity for $|\xi|\ge 1$ we find
\begin{eqnarray}
   &&\int_{|\xi|\ge|\mu|} a(x,\xi)\, \log^t|\xi|\, b(x,\xi,\mu)d\xi\nonumber\\
   &=&|\mu|^{a+n+j-mN} \int_{|\xi|\ge 1} a(x,\xi)\, (\log |\xi\mu|)^t\,
       b(x,\xi,\frac{\mu}{|\mu|})d\xi\\
  &=&\sum_{\sigma=0}^t  {t\choose \sigma} \int_{|\xi|\ge 1} a(x,\xi)
   \,\log^{t-\sigma}|\xi| \, b(x,\xi,\frac{\mu}{|\mu|})\,d\xi\:
   |\mu|^{a+n-j-mN}\,\log^\sigma|\mu|.\nonumber
\end{eqnarray}
Since $b$ is holomorphic in $\mu$ by \cite[Lemma 2.3]{\GS} this is
actually an expansion in terms of the functions
$\mu^{a+n-j-mN} \,\log^\sigma \mu$.

Now we expand using \cite[Theorem 1.12]{\GS}
\begin{equation}
  b(x,\xi,\mu)=\sum_{0\le \nu<M} b_\nu(x,\xi) \mu^{-m\nu-Nm}+R_M(x,\xi,\mu)
\end{equation}
where $b_\nu$ is homogeneous in $|\xi|$ for $|\xi|\ge 1$ and
\begin{equation}
  R_M(x,\xi,\mu)=\mbox{\large O}\big(<\xi>^{mM-l'}\mu^{-m(M+N)}\big),
\end{equation}
thus
\begin{eqnarray}\begin{array}{rl}
   &\DST\int_{|\xi|\le 1} a(x,\xi)\, \log^t|\xi|\, b(x,\xi,\mu)d\xi\\[1em]
   =&\DST\sum_{0\le \sigma<M} \mu^{-m\nu-Nm}
      \int_{|\xi|\le 1} a(x,\xi)\, \log^t |\xi| \,b_\nu(x,\xi) d\xi
     +O(\mu^{-m(M+N)}).\end{array}
\end{eqnarray}
Furthermore,
\begin{eqnarray}
   &&\int_{1\le |\xi|\le |\mu|} a(x,\xi)\, \log^t|\xi|\, b_\nu(x,\xi)d\xi\,
     \mu^{-m\nu-Nm}\nonumber\\
  &=&\mu^{-m\nu-Nm} c_\nu(x) \int_1^{|\mu|} r^{a-l-|\ga|+m\nu-l'+n-1}
     \log^t r dr\nonumber \\
  &=&\left\{\begin{array}{ll}
        \DST \sum_{\sigma=0}^t c_{\nu,\sigma}(x) (\log^\sigma|\mu|)
     \mu^{-m\nu-Nm}|\mu|^{a-j+m\nu+n},&\DST a-j+m\nu+n\not=0,\\
      \DST \frac{c_\nu(x)}{t+1} (\log^{t+1}|\mu|) \mu^{-m\nu-Nm},&
      \DST a+m\nu+n-j=0.
   	   \end{array}\right.\mylabel{G3-4.3}
\end{eqnarray}
Again invoking \cite[Lemma 2.3]{\GS} we end up with the desired expansion.
The remainder term coming from $R_M$ is estimated exactly as in
\cite[p. 498]{\GS}. \myref{G3-4.3} also shows that
$c_{j,k+1}=0$ if $\frac{j-a-n}{m}\not\in \Z_+$.\endproof   



As a consequence of Theorem \plref{S1-3.5} we have an asymptotic
expansion 
\begin{equation}
   \Tr(A(P-\gl)^{-N}))\sim \sum_{j=0}^\infty \sum_{l=0}^{k+1}
       c_{jl} \gl^{\frac{n+a-j}{m}-N} \log^l\gl
     \quad + \sum_{j=0}^\infty d_j \gl^{-j-N},
   \mylabel{G1-3.8}
\end{equation}
where $c_{j,k+1}=0$ if $\frac{j-a-n}{m}\not\in\Z_+$.

If the eigenvalues of the leading symbol of $P$ lie in $\Re \gl>0$ then
via the appropriate Cauchy integral we obtain the heat expansion
\begin{equation}
   \Tr(A e^{-tP})\sim_{t\to 0+} \sum_{j=0}^\infty t^{\frac{j-n-a}{m}}
     \tilde c_j(\log t)+ \sum_{j=0}^\infty \tilde d_j t^j 
   \mylabel{G1-3.9}
\end{equation}
with a polynomial $\tilde c_j\in \C[x]$ of degree
\begin{equation}
   \deg \tilde c_j\le \casetwo{k}{\frac{j-a-n}{m}\not\in\Z_+,}%
     {k+1}{\frac{j-a-n}{m}\in\Z_+.}
\end{equation}

Furthermore, the generalized $\zeta$--function
\begin{equation}
   \Tr(A P^{-s})=\frac{1}{\Gamma(s)}\int_0^\infty t^{s-1}\Tr(Ae^{-tP})dt
   \mylabel{G1-3.10}
\end{equation}
has a meromorphic continuation to $\C$ with poles in
$\{ \frac {a+n-j}{m}\,|\, j\in \Z_+\}$
of order k+1. Cf. e.g. \cite[Lemma 2.1]{BruLes:SGAC}, \cite{GruSee:ZEFAPSO},
\cite[Sec 2.1]{Les:OFTCSAM}, for a discussion of the relation
between the expansions \myref{G1-3.8}, \myref{G1-3.9} and the poles
of \myref{G1-3.10}.

Using the notation $f(s)=:\DST \sum_k \TST \Res_k f(s_0) (s-s_0)^{-k}$
for Laurent expansions we note explicitly
\begin{eqnarray}
     \Res_{k+1} \Tr(AP^{-s})|_{s=\frac{a+n-j}{m}}&=&
        \frac{(-1)^k}{\Gamma(\frac{a+n-j}{m})}\frac{d^k}{dt^k}\tilde 
   c_j(t)|_{t=0},
       \; \TST\frac{j-n-a}{m}\not\in\Z_+,\mylabel{G1-3.11}\\
      \Res_{k+1} \Tr(AP^{-s})|_{s=-j}&=& (-1)^{k+j+1} j!
      \frac{d^{k+1}}{dt^{k+1}}\tilde c_{a+n+mj}(t)|_{t=0},
\;\TST  \frac{j-n-a}{m}\in\Z_+.\mylabel{G1-3.12}
\end{eqnarray}

\section[The \poly\ noncommutative residue]{The \poly\ noncommutative\\ residue}

In this section we consider a compact closed manifold $M$ of dimension $n$.
Let $E$ be a $C^\infty$ vector bundle over $M$ and $A\in \CL^{a,k}(M,E)$ a
pseudodifferential operator with \poly\ symbol.

We pick an elliptic classical \pdo\ \op\ $P\in \CL^{m,0}(M,E)$ of order $m>0$
whose leading symbol is {\em scalar} and {\em positive}. For instance
we could take a generalized Laplace operator for $P$. 
The assumption that the leading symbol of $P$ is scalar guarantees
that the commutator $[P,A]$ lies in $\CL^{a+m-1,k}(M,E)$, which makes
life easier in the sequel. We emphasize however that this assumption is
not really important.

We put
\begin{equation}
        \nabla_P^0 A:=A,\quad \nabla_P^{j+1}A:=[P,\nabla^j A],
     \mylabel{G1-4.1}
\end{equation}
and by induction we have
\begin{equation}
	\nabla_P^j A\in \CL^{a+j(m-1),k}(M,E).
     \mylabel{G1-4.2}	
\end{equation}
\begin{dfn}{S1-4.1} We put
\begin{eqnarray*}
      \Res_k(A,P)&:=& m^{k+1} \Res_{k+1} \Tr(AP^{-s})|_{s=0}\\
     &=& m^{k+1} (-1)^{k+1} \tilde c^{(k+1)}_{a+n}(0)\\
    &=& m^{k+1}(-1)^{k+1} (k+1)!\times\;\mbox{coefficient of }\;
      \log^{k+1} t\;\mbox{in the asymptotic}\\
    &&\mbox{expansion of}\;\Tr(Ae^{-tP})\;\mbox{as}\; t\to 0.\\
    &=& m^{k+1}(-1)^N (k+1)!\times\;\mbox{coefficient of }\;
      \gl^{-N}\log^{k+1} \gl\;\mbox{in the asymptotic}\\
    &&\mbox{expansion of}\;\Tr(A(P-\gl)^{-N})\;\mbox{as}\; \gl\to \infty.
\end{eqnarray*}
\end{dfn}

We abbreviate
\begin{equation}
   \Delta_N:=\{ (t_1,\ldots,t_N)\in \R^N\,|\, 0\le t_1\le \ldots \le
     t_N\le 1\}.
   \mylabel{G1-4.3}
\end{equation}
Then we have the well--known formula\footnote{Actually it was not well--known
to the author; he has learned the formula and the following lemma from
Henri Moscovici}
\begin{equation}
  e^{-tP}A =\sum_{j=0}^{N-1} \frac{(-t)^j}{j!} (\nabla_P^j A) 
       e^{-tP}+R_N(A,P,t),
    \mylabel{G1-4.4}
\end{equation}
where
\begin{equation}\begin{array}{rcl}
  \DST R_N(A,P,t)&=&\DST (-1)^N\int_{t\Delta_N}
       e^{-t_1P} (\nabla_P^N A) e^{-(t-t_1)P}\dtn{N}\\[1.5em]
   &=&\DST\frac{(-t)^N}{(N-1)!}\int_0^1 (1-s)^{N-1}
       e^{-stP} (\nabla_P^N A) e^{-(1-s)tP}ds.
   		\end{array}
    \mylabel{G1-4.5}
\end{equation}
\begin{lemma}{S1-4.2} If $p,\eps,N>0$ such that $(N-a)/m-p-\eps>0$ then
we have the estimate
$$ \|R_N(A,P,t)(P+c)^p\|\le C t^{\frac{N-a}{m}-p-\eps},\quad
0<t\le 1.$$
Here $c>0$ is any constant such that $P+c$ is invertible.
\end{lemma}
\proof Since the leading symbol of $P$ is positive the operator
$P+c$ is invertible for some $c>0$. 
$\nabla_P^N A\in \CL^{a+N(m-1),k}(M,E)\subset \rmL^{a+N(m-1)+m\eps}(M,E)$, hence
\newline 
$(P+c)^{\frac {N-a}{m}-N-\eps}\nabla_P^N A$ is bounded and thus we obtain the
estimate
\begin{eqnarray*}
   \|e^{-t_1 P} (\nabla_P^N A )  e^{-(t-t_1)P}(P+c)^p\| &\le&
      C \| (P+c)^{\frac {a-N}{m}+N+p+\eps} e^{-(t-t_1)P}\|\\
    &\le& C (t-t_1)^{\frac {N-a}{m}-N-p-\eps},\quad 0\le t_1< t \le 1.
\end{eqnarray*}
Integrating this inequality gives the desired estimate.\endproof

\begin{prop}{S1-4.3} Let $A\in \CL^{a,k}(M,E)$ and let
$P\in \CL^{m,0}(M,E)$ be a classical elliptic \pdo\ \op\ of order
$m>0$ whose leading symbol is scalar and positive.
Then 
\begin{enumerate}
\renewcommand{\labelenumi}{{\rm \arabic{enumi}.}}
\item $\Res_k(A,P^\ga)=\Res_k(A,P)$ for any $\ga>0$.
\item If $P_u$ is a smooth 1--parameter family of such operators
then $\Res_k(A,P_u)$ is independent of the parameter $u$.
\item If $B\in \CL^{b,l}(M,E)$ then $\Res_{k+l}([A,B],P)=0$.
\end{enumerate}
\end{prop}
\proof 1. is straightforward.

2. In view of Lemma \plref{S1-4.2} we have for $N$ large enough
\begin{equation}\begin{array}{rcl}
 \DST\frac{d}{du} A e^{-tP_u}&=&\DST-\int_0^t A e^{-(t-t_1) P_u}
        (\frac{d}{du} P_u) e^{-t_1 P_u}dt_1\\[1.5em]
     &=&\DST \sum_{j=0}^{N-1} \frac{(-t)^{j+1}}{(j+1)!} (A \nabla_P^j 
     \frac{d}{du} P_u)e^{-tP}+R_N(t),
   		\end{array}
\end{equation}
where
\begin{equation}
\|R_N(t)\|_{\rm tr} = O(t),\quad t\to 0.
\end{equation}
Here $\|\cdot\|_{\rm tr}$ denotes the trace norm and this trace norm
estimate follows from Lemma \plref{S1-4.2} since $(P+c)^{-p}$ is
trace class if $p>\dim M/m$.

Thus we conclude
\begin{equation}
 \frac{d}{du} \Tr(A e^{-tP_u})=
   \sum_{j=0}^{N-1} \frac{(-1)^{j+1}}{(j+1)!}t^{j+1}\Tr(A (\nabla_P^j 
     \frac{d}{du} P_u)e^{-tP})+O(t),\quad t\to 0.
\end{equation}
Since $\Tr(A \nabla_P^j (\frac{d}{du} P_u)e^{-tP})$ has no
$t^{-j-1}\log^{k+1} t$ term in its asymptotic expansion (cf. \myref{G1-3.9})
we reach the conclusion.

3. We use again Lemma \plref{S1-4.2} and we find for $N$ large enough
\begin{equation} Ae^{-tP}B= 
    \sum_{j=0}^{N-1} \frac{(-t)^j}{j!} A (\nabla_P^j B) 
      e^{-tP}+R_N(t),
\end{equation}
where $\|R_N(t)\|_{\rm tr} = O(t),\; t\to 0.$ Thus
\begin{equation}
\Tr([A,B]e^{-tP})=- \sum_{j=1}^{N-1} \frac{(-1)^j}{j!} t^j 
   \Tr(A (\nabla_P^j B) e^{-tP})+O(t),\quad t\to 0.
\end{equation}
Since $A \nabla_P^j B\in \CL^{a+b+j(m-1),k+l}(M,E)$ the highest
$\log t$ power occuring as coefficient of $t^{-j}$ in the asymptotic
expansion of  $\Tr(A (\nabla_P^j B) e^{-tP})$ as $t\to 0$ is
$\le k+l$ and hence $\Res_{k+l}([A,B],P)=0$.\endproof

This proposition immediately implies:

\begin{theorem}{S1-4.4} For $A\in \CL^{a,k}(M,E)$ choose $P\in \CL^{m,0}(M,E)$
as in the preceding proposition. Then
$$\Res_k(A):=\Res_k(A,P)$$
is well--defined independent of the particular $P$ chosen. $\Res_k$
is a linear functional on $\CL^{*,k}(M,E)$ which vanishes on
appropriate commutators. More precisely,
\begin{equation}
   \Res_{k+l}([A,B])=0
   \mylabel{G1-4.6}
\end{equation}
for $A\in \CL^{a,k}(M,E), B\in \CL^{b,l}(M,E)$.
\end{theorem}
\proof Given $P_j\in\CL^{m_j,0}(M,E)$, $j=1,2$, in view of 1. of the preceding
proposition we may replace $P_j$ by $P_j^{1/m_j}$ to obtain operators of
order 1. Then $P_u:=uP_1+(1-u)P_2$, $0\le u\le 1$,
is a smooth 1--parameter family
to which 2. of the preceding proposition applies. Hence
$\Res_k(A,P_1)=\Res_k(A,P_2).$ \myref{G1-4.6} now follows from
3. of loc. cit.\endproof

We note that $\Res_0$ coincides with the noncommutative residue
of {\sc Wodzicki}. This follows directly from the definition and
\cite[Thm. 1.4]{Kas:RNC}.

Although \myref{G1-4.6} holds $\Res_k$ is not a trace on the full algebra
$\CL^{*,*}(M,E)$. We will see in Theorem \plref{S1-4.9} below 
that this algebra does not have any nontrivial traces.
However, as we will see at the end of this section the sequence of
functionals $(\Res_k)_{k\in\Z_+}$ corresponds to a trace
on a graded algebra constructed from the filtration
$(\CL^{*,k})_{k\in\Z_+}$. 

Note furthermore that $\Res_k(A)=0$ if ${\rm ord}(A)<-\dim M$ because then
$A$ is of trace class and hence $\Tr(AP^{-s})$ is regular at $s=0$.

\subsection{The residue density}

\begin{prop}{S1-4.5}{\rm (cf. \cite[Prop. 3.2]{Wod:NRCIF})}
Let $M$ be a (not necessarily compact) smooth manifold.
For $A\in \CL^{m,k}(M,E,F)$ there exists a density
$$\go_k(A)\in\cinf{M,\Hom(E,F)\otimes|\Omega|}$$
with the following properties:
\begin{enumerate}
\renewcommand{\labelenumi}{{\rm \arabic{enumi}.}}
\item $A\mapsto \go_k(A)$ is $\C$--linear.
\item If $\varphi\in\cinfz{M}$ then $\go_k(\varphi A)=\varphi \go_k(A)$.
\item If $\kappa:U\to V\subset \R^n$ is a local chart then
$$\kappa^*(\go_k(\kappa_*A))=\go_k(A).$$
\item In a local coordinate chart we have
\begin{eqnarray}
\go_k(A)(x)&=& \frac{(k+1)!}{(2\pi)^n}\res_k(a_{-n}(x,\cdot))|dx|\nonumber \\
    &=& \frac{(k+1)!}{(2\pi)^n} \Big(\int_{|\xi|=1} a_{-n,k}(x,\xi) |d\xi|
\Big) |dx|.
        \mylabel{G1-4.7}
\end{eqnarray}
\item If the complete symbol of $A\in \CL^{m,k}(M,E,F)$ vanishes at
$p\in M$ then $\go_k(A)(p)=0$. If $m<-\dim M$ then $\go_k(A)=0$.
\item Let $E=F$, $A\in \CL^{a,k}(M,E), B\in \CL^{b,l}(M,E)$.
If $A$ or $B$ have compact support then
$$\int_M \tr_{E_x}(\go_{k+l}([A,B])(x))=0.$$
\end{enumerate}
\end{prop}
\proof We take \myref{G1-4.7} as the definition of $\go_k(A)$. 
To show that $\go_k(A)$ is well defined we proceed along the
lines of \cite[Thm. 1.4]{\Schretal}.

Let $U\subset \R^n$, $A\in \CL^{m,k}(U,E,F)$ and $\kappa:U\to V$ a
diffeomorphism. Since this consideration is local we may assume $E,F$
to be trivial bundles. By Proposition \plref{S1-3.4} we have
$$ \sigma_{\gk_*A}(\kappa(x),\eta)\sim
     \sum_{\ga\in \Z_+^n}  (\pl_\xi^\ga \sigma_A)(x,
      D\gk^t(x) \eta) \Phi_\ga(x,\eta),$$
where $\Phi_\ga(x,\eta)$ is a polynomial in $\eta$ of degree not higher
than $|\ga|/2$ and $\Phi_0=1$. Hence in view of \myref{res-G1.10},
\myref{res-G1.11}, and Lemma \plref{res-S1.8}
\begin{eqnarray*}
   && \res_k[ \sigma_{\gk_*A}(\kappa(x),\eta)_{-n}] |d\gk(x)|\\
   &=& |\det D\gk(x)| \sum_{\ga\in\Z_+^n} 
    \res_k\big[\big((\pl_\xi^\ga \sigma_A)(x,
      D\gk^t(x) \eta) \Phi_\ga(x,\eta)\big)_{-n}\big]|dx|\\
  &=&  |\det D\gk(x)| \res_k[\sigma_A(x,D\gk^t(x) \eta)_{-n}]|dx|\\
  &=& \res_k [ \sigma_A(x,\eta)_{-n}] |dx|.
\end{eqnarray*}

This proves 3. and 4. Furthermore, 1., 2., and 5. are simple consequences
of 4. 

To prove 6., using a partition of unity we may assume that $M=U$ is a
coordinate chart and $A$ has compact support.
Now we proceed exactly as in the proof of 
\cite[Thm. 1.4]{\Schretal}: denote by $a,b$ the complete symbols of $A,B$.
We obtain from the symbol expansion of a product
Proposition \plref{S1-3.3}
\begin{eqnarray*}
  &&  \frac{(2\pi)^n}{(k+1)!}\tr_{E_x}\go_{k+l} ([A,B])(x)\\
&=&   \int_{|\xi|=1} \Big(\sum_{\ga\in\Z_+^n}\frac{i^{-|\ga|}}{\ga!} 
    \tr_{E_x}(\pl_\xi^\ga a\,\pl_x^\ga b-\pl_\xi^\ga b\,
   \pl_x^\ga a)\Big)_{-n,k}
    |d\xi| |dx|\\
   &=& \int_{|\xi|=1} \Big(\sum_{\ga\in\Z_+^n}\frac{i^{-|\ga|}}{\ga!} 
    \tr_{E_x}(\pl_\xi^\ga a\,\pl_x^\ga b-\pl_x^\ga a\,
   \pl_\xi^\ga b)\Big)_{-n,k}
    |d\xi| |dx| .
\end{eqnarray*}
Now, as noted in loc. cit.
\begin{equation}
     \pl_\xi^\ga a\,\pl_x^\ga b-\pl_x^\ga a\,\pl_\xi^\ga b
\end{equation}
can be written as a sum of derivatives
\begin{equation}
    \sum_{j=1}^n \pl_{\xi_j} A_j+\pl_{x_j} B_j,
\end{equation}
where $A_j, B_j$ are bilinear expressions in $a, b$ and their derivatives
and $a$ or one of its derivatives occurs in every summand. Hence, since
$a$ has compact support the $A_j, B_j$ have compact support. The assertion
now follows from Lemma \plref{res-S1.8}.\endproof

\begin{lemma}{S1-4.6} Let $M$ be a compact Riemannian manifold, $\dim M=n$.
Furthermore, let $Q\in \CL^{-n,k}(M,E)$ with leading symbol
$q_n(x,\xi)=|\xi|^{-n} \log^k|\xi|$. Then
$$ \Res_k(Q)=\int_M \go_k(Q)\not=0.$$
\end{lemma}
\proof This lemma could be derived from the proof of
\cite[Thm 2.1]{\GS}. But since our situation is much simpler
we will give an ad hoc proof.

Note first that in a coordinate system
$$\go_k(Q)(x)=\frac{(k+1)!}{(2\pi)^n}\vol(S^{n-1})d\vol_M(x)$$
and thus
$$\int_M \go_k(Q)=\frac{(k+1)!}{(2\pi)^n}\vol(S^{n-1})\vol(M)\not=0.$$

We pick an elliptic operator $P\in \CL^{m,0}(M,E), m>0$, whose leading
symbol $p(x,\xi)$ is scalar and positive. 

Then since $Q$ has order $-n$ the coefficient of $\gl^{-N}\log^{k+1}\gl$
in the asymptotic expansion of $\Tr(A(P-\gl)^{-N})$ equals
the corresponding coefficient of
$$(2\pi)^{-n}\int_{T^*M} \varphi(|\xi|) |\xi|^{-n} \log^k|\xi|
   (p(x,\xi)-\gl)^{-N} d\xi dx.$$

Here, $\varphi\in\cinf{\R}$ is a function with $\varphi(t)=0$, if
$t\le 1/2$ and $\varphi(t)=1$, if $t\ge 1$.

Now
\begin{eqnarray*}
    &&\int_{\R^n} \varphi(|\xi|)|\xi|^{-n}\log^k|\xi|(p(x,\xi)-\gl)^{-N}d\xi\\
   &=&\gl^{-N} \int_{\R^n} 
       \varphi(|\gl^{1/m}\xi|)|\xi|^{-n}\log^k|\gl^{1/m}\xi|
      (p(x,\xi)-1)^{-N}d\xi\\
  &=& \gl^{-N} \int_{|\xi|\le 1} 
       \varphi(|\gl^{1/m}\xi|)|\xi|^{-n}\log^k|\gl^{1/m}\xi|
      (p(x,\xi)-1)^{-N}d\xi+O(\gl^{-N} \log^k \gl)\\
  &=&  \gl^{-N} \int_0^1\int_{|\xi|= 1} 
       \varphi(|\gl^{1/m}\varrho|)\varrho^{-1}\log^k|\gl^{1/m}\varrho|
      (p(x,\varrho\xi)-1)^{-N}d\xi d\varrho+O(\gl^{-N} \log^k \gl)\\
  &=& (-1)^N \gl^{-N} (\log^k \gl^{1/m}) \vol(S^{n-1})
      \int_0^1 \varphi(|\gl^{1/m}\varrho|)\varrho^{-1}d\varrho
   +O(\gl^{-N} \log^k \gl)\\
  &=& \frac{(-1)^N}{m^{k+1}}\vol(S^{n-1}) \gl^{-N} \log^{k+1} \gl
      +O(\gl^{-N} \log^k \gl)
\end{eqnarray*}
and hence
\begin{eqnarray*}
&&(2\pi)^{-n}\int_{T^*M} \varphi(|\xi|) |\xi|^{-n} \log^k|\xi|
   (p(x,\xi)-\gl)^{-N} d\xi dx\\
&=& \Big(\frac{(-1)^N}{m^{k+1}(k+1)!}
    \int_M\go_k(Q)\Big) \gl^{-N} \log^{k+1} \gl
      +O(\gl^{-N} \log^k \gl)
\end{eqnarray*}
which proves the assertion in view of Definition \plref{S1-4.1}.\endproof

\begin{prop}{S1-4.7} Let $M$ be a connected compact manifold of
dimension $n>1$.
Choose a $Q\in \CL^{-n,k}(M)$ with $ \Res_k(Q)=\int_M \go_k(Q)\not=0$.
Then there exist $P_1,\ldots,P_N\in \CL^{1,0}(M)$ such that
for any $A\in \CL^{a,k}(M)$ there exist 
$Q_1,\ldots,Q_N\in \CL^{a,k}(M)$ such that
$$A-\sum_{j=1}^N [P_j,Q_j]- \frac{\int_M \go_k(A)}{\Res_k(Q)} Q\in \rmL^{-\infty}(M).$$
\end{prop}

\proof We choose $P_1,\ldots,P_N\in \CL^{1,0}(M)$ such that
the differentials of the leading symbols span $T^*M\setminus 0$ at every point.
We consider the leading symbol $\sigma_\rmL^a(A)\in\cp^{a,k}(T^*M)$ of
$A$. If $a\not=-n$ then by Theorem \plref{res-S2.5} there exist
$Q_1^{(1)},\ldots,Q_N^{(1)}\in \CL^{a,k}(M)$ such that
   $$A-\sum_{j=1}^N [P_j,Q_j^{(1)}]\in \CL^{a-1,k}(M).$$
(Note that the leading symbol of a commutator is the Poisson bracket
of the leading symbols).

If $a-1\not=-n$ we iterate the procedure. Thus if $a\not\in\{l\in \Z\,|\,
 l\ge -n\}$ then  by induction we find operators $Q_j^{(l)}\in \CL^{a,k}(M)$
such that
\begin{equation}
A^{(l)}:=A-\sum_{j=1}^N [P_j,Q_j^{(l)}]\in \CL^{a-l,k}(M).
  \mylabel{G1-4.8}
\end{equation}
If $a\in\{l\in \Z\,|\, l\ge -n\}$ then \myref{G1-4.8} holds for $l\le a+n$.
In view of Proposition \plref{S1-4.5} we have
$$\int_M \go_k(A^{(a+n)}-\frac{\int_M \go_k(A)}{\Res_k(Q)} Q)=0.$$
Using again Theorem \plref{res-S2.5} we find by induction
operators  $Q_j^{(l)}\in \CL^{a,k}(M)$
such that
\begin{equation}
A^{(l)}:=A-\sum_{j=1}^N [P_j,Q_j^{(l)}]-\frac{\int_M \go_k(A)}{\Res_k(Q)} Q
  \in \CL^{a-l,k}(M).
  \mylabel{G1-4.9}
\end{equation}
Now choose $Q_j\in \CL^{a,k}(M)$ with $Q_j-Q_j^{(l)}\in \CL^{a-l,k}(M)$.
Then we reach the conclusion. 
\endproof

\begin{cor}{S1-4.8} For $A\in \CL^{a,k}(M,E)$ we have
$$   \Res_k(A)=\int_M \tr_{E_x}\go_k(A)(x)= 
   \frac{(k+1)!}{(2\pi)^n} \int_{S^*M} a_{-n,k}(x,\xi) |d\xi dx|.$$
\end{cor}
\proof For $E=\C$ 
this is an immediate consequence of Theorem \plref{S1-4.4},
Proposition \plref{S1-4.5} and the preceding proposition.
For arbitrary $E$ it follows from the facts that there exists
a bundle $F$ making $E\oplus F$ trivial and 
$\CL^{a,k}(M,\C^r)\cong \CL^{a,k}(M)\otimes M_r(\C)$.
\endproof

We note explicitly as a consequence of Proposition \plref{S1-4.7} that
if $A\in \CL^{a,k}(M)$ there are always $Q_j\in \CL^{a,k+1}(M)$ such 
that
$$A-\sum_{j=1}^N [P_j,Q_j]\in \rmL^{-\infty}(M).$$

This can actually be improved, using the fact that every classical
pseudodifferential operator with vanishing residue is a sum
of commutators. This result is due to
{\sc Wodzicki} {\rm \cite{Wod:SALI}}, a generalization to algebras of
Fourier integral operators
is due to {\sc Guillemin} \cite{Gui:RTCAFIO}.
In particular every smoothing operator
is a sum of commutators of classical pseudodifferential operators. 
Since \cite{Wod:SALI} is written
in Russian we briefly sketch the argument: one actually has to show
that $\CL^*(M)/[\CL^*(M),\CL^*(M)]$ is one--dimensional. Given any trace
$\tau$ on $\CL^*(M)$ consider first $\tau|\rmL^{-\infty}(M)$. By
\cite[Appendix]{Gui:RTCAFIO} 
$\rmL^{-\infty}(M)/[\rmL^{-\infty}(M),\rmL^{-\infty}(M)]\cong \C$ 
and it is spanned by the $\rmL^2$--trace. Thus 
$\tau|\rmL^{-\infty}(M)=
c\, \Tr_{\rmL^2}$. The same argument as in the introduction
\myref{G1-I.1} then shows $c=0$. Hence $\tau$ induces a trace on
$\CL^*(M)/\rmL^{-\infty}(M)$ and thus in view of Proposition \plref{S1-4.7}
it is a constant multiple of the {\sc Wodzicki} residue.

Summing up we can improve Proposition \plref{S1-4.7} as follows:

\begin{prop}{S2-4.9} Let $M$ be a connected compact manifold of
dimension $n>1$.
Choose a $Q\in \CL^{-n,k}(M)$ with $ \Res_k(Q)=\int_M \go_k(Q)\not=0$.
Then there exist $P_1,\ldots,P_N\in \CL^{1,0}(M)$ such that
for any $A\in \CL^{a,k}(M)$ there exist 
$Q_1,\ldots,Q_N\in \CL^{a,k}(M)$ and classical pseudodifferential operators
$R_1,\ldots,R_N, S_1,\ldots, S_N\in \CL^{*,0}(M)$ such that
$$A=\sum_{j=1}^N [P_j,Q_j]+\sum_{j=1}^N [R_j,S_j]+
     \frac{\int_M \go_k(A)}{\Res_k(Q)} Q.$$
In particular, since $\Res_{k+1}(A)=0$ there exist 
$Q_1,\ldots,Q_N\in \CL^{a,k+1}(M)$ and \\
$R_1,\ldots,R_N,$ $S_1,\ldots, S_N \in \CL^{*,0}(M)$ such that
$$A=\sum_{j=1}^N [P_j,Q_j]+\sum_{j=1}^N [R_j,S_j].$$
\end{prop}

\begin{remark}
This result generalizes the corresponding result for classical
pseudodifferential operators due to {\sc Wodzicki} \cite{Wod:SALI}
\cite[Prop. 5.4]{Kas:RNC}. In \cite{Gui:RTCAFIO}, {\sc Wodzicki}'s result
was generalized to certain algebras of Fourier integral operators.
\end{remark}

An immediate consequence is the 

\begin{cor}{S1-4.9} There are no nontrivial traces
on the algebra $\CL^{*,*}(M)$. 
\end{cor}
\comment{\proof Let $\tau$ be a trace on $\CL^{*,*}(M,E)$.
We pick a bundle $F$ such that $E\oplus F=\C^r_M$,
the trivial bundle of rank $r$. Then $\tilde \tau:=\tau\oplus 0$ is
a trace on $\CL^{*,*}(M,\C^r_M)=\CL^{*,*}(M)\otimes M_r(\C)$.
Then it is easy to check that
$$\tilde \tau=\bar \tau \otimes \tr,$$
where $\tr$ is the trace on $M_r(\C)$ and $\bar \tau$ is a trace
on $\CL^{*,*}(M)$ with the same properties as $\tau$. 
Hence it suffices to show that
$\bar \tau=0$. But this is clear since every element of $\CL^{*,*}(M,E)$
is a sum of commutators.\endproof}

\comment{Now pick a surjective Fredholm operator
$T\in \CL^{0,0}(M,\C^r)$ of index one, which exists by the
Atiyah--Singer index theorem. Then $TT^*=I$ and
$I-T^*T\in \rmL^{-\infty}(M,\C^r)$ is a rank one projection. Thus
$$\tau\otimes \tr_{\C^r}(I-T^*T)=0$$ and
hence $c=0$. }

\comment{Next consider $A\in \CL^{*,*}(M)$. In view of Proposition
\plref{S1-4.7} and the preceding remark there exist $P_1,\ldots, P_N\in \CL^{1,0}(M)$
and $Q_1,\ldots, Q_N\in \CL^{*,*}(M)$ such that
$$A-\sum_{j=1}^N [P_j,Q_j]\in \rmL^{-\infty}(M)$$
and thus
$$\tau(A)=\tau(A-\sum_{j=1}^N [P_j,Q_j])=0.\epformel$$}

\comment{The following lemma which was used in the previous proof is probably
well--known:

\begin{lemma}{S1-4.10}  Let $\tau$ be a trace on the algebra
$\rmL^{-\infty}(M,E)$. If $\tau|\rmL^{-\infty}(M,E)$ is continuous with
respect to the $C^\infty$--topology on $\rmL^{-\infty}(M,E)\cong
  \cinf{M\times M, E\boxtimes E^*}$ then $\tau=c\, \tr_{\rmL^2}$
for some constant $c$.
\end{lemma}
\proof We choose a Riemannian metric and equip the bundle $E$ with
a hermitian structure. Since every operator $T\in \rmL^{-\infty}(M,E)$
can be approximated in the $C^\infty$--topology by operators of
finite rank it suffices the prove $\tau=c\, \tr_{\rmL^2}$ on
finite rank operators  
in $\rmL^{-\infty}(M,E)$. By linearity it suffices to consider rank one
operators. But this is
an exercise in linear algebra:
introduce the rank one operators $T_{f,g}(u):= (g,u) f$,
$f,g\in \cinf{E}$. Then
$$T_{f,g}\circ T_{f_1,g_1}=(g,f_1) T_{f,g_1}$$
and $T_{g,g}$ is   a projection if $\|g\|=1$. Thus if $\|f\|=\|g\|=1$
then
$$T_{f,g}\circ T_{g,f}=T_{f,f},\quad T_{g,f}\circ T_{f,g}=T_{g,g},$$
and thus $\tau$ is constant on rank one projections.
Let $c$ be the common value. For arbitrary $f,g\not=0$ we then find
\begin{epeqnarray}
\tau(T_{f,g})&=& \tau(\frac{1}{\|g\|^2} T_{f,g}\circ T_{g,g})
  =\frac{1}{\|g\|^2} \tau(T_{g,g}\circ T_{f,g})\\
 &=&(g,f) \frac{1}{\|g\|^2} \tau(T_{g,g})=c(g,f)=c \,\tr_{\rmL^2}(T_{f,g}).
\end{epeqnarray}}

We close this section with two further remarks. First we note
that our algebra $\CL^{*,*}(M,E)$ provides examples of spectral
triples with discrete dimension spectrum 
of infinite multiplicity as defined by
{\sc Connes} and {\sc Moscovici} \cite[Def. II.1]{ConMos:LITNG}.
The algebra involved in a spectral triple consists of
bounded operators. Therefore, we put
\begin{equation}
    \ca(M,E):=\{A\in\CL^{0,*}(M,E)\,|\, \sigma_L(A)\in\cp^{0,0}\}.
\end{equation}
This is a subalgebra of $\CL^{*,*}$ consisting of bounded operators.

If $P\in \CL^{1,0}(M,E)$ is a classical elliptic
pseudodifferential operator whose leading symbol is scalar and
positive then \myref{G1-3.10} shows that 
\begin{equation}
     (\ca(M,E),\rmL^2(M,E),P)
\end{equation}
is a spectral triple with dimension spectrum
${\rm Sd}:=\{k\in\Z\,|\, k\le \dim M\}$ of infinite multiplicity.

The assumption that the leading symbol is scalar does not apply
to the Dirac operator. However, this assumption 
is needed only to guarantee
that $[P,A]$ is of order $0$. Thus, given a spin manifold, $M$,
with Dirac operator, $D$, acting on $\cinf{M,S}$
we denote by $\ca_{\rm scal}(M,S)$ the algebra of operators
$A\in\ca(M,S)$ whose complete symbol is scalar. Having a scalar
complete symbol is a coordinate invariant property. Hence
\begin{equation}
     (\ca_{\rm scal}(M,S),\rmL^2(M,S),D)
\end{equation}
is another spectral triple with dimension spectrum 
${\rm Sd}:=\{k\in\Z\,|\,k\le\dim M\}$.

Our second remark concerns the fact that 
although $\Res_{k+l}([A,B])=0$
if $A\in\CL^{a,k}$, $B\in \CL^{b,l}$, the functional $\Res_k$ is not
a trace on the whole algebra $\CL^{*,*}$. To shed some light on this fact
we consider an arbitrary filtered algebra 
\begin{equation}
     \cb:=\bigcup_{k=0}^\infty \cb^k, \quad \cb^k\subset \cb^{k+1}.
     \label{G5-5.1}
\end{equation}
From $\cb$ we can construct a graded algebra
\begin{equation}
     \cg\cb:=\moplus_{k=0}^\infty \cb^k,
\end{equation}
where the product is given by
\begin{equation}
    ((a_k)_{k\ge 0}\circ (b_k)_{k\ge 0})_m:=\sum_{j=0}^m a_j b_{m-j}.
\end{equation}
Then it is straightforward to see that traces on $\cg\cb$ are in
one--one correspondence to sequences of linear functionals
$\tau_k:\cb^k\longrightarrow \C$ satisfying
\begin{equation}
     \tau_{k+l}([A,B])=0
\end{equation}
for $A\in\cb^k, B\in \cb^l$.

Thus Proposition \plref{S2-4.9} says that if $M$ is connected, $\dim M>1$,
then up to a scalar factor there is exactly one trace on the algebra
$\cg\CL^{*,*}(M)$ constructed from the filtration 
\begin{equation}
        \CL^{*,*}(M)=\bigcup_{k=0} ^\infty \CL^{*,k}(M).
\end{equation}
This fact was communicated to the author by {\sc R. Nest.}

\section{The Kontsevich--Vishik trace}

\newcommand{\TR}{{\rm TR}}

In this section we briefly show that the analogue of the {\sc Kontsevich--Vishik}
\cite{KonVis:DEPO,KonVis:GDEO}
trace exists on $\CL^{*,*}$, too. During the whole section
$M$ will be a compact manifold without boundary.
We will present two proofs which
slightly differ from the method of loc. cit. 

The first proof is exactly along the lines of Definition \plref{S1-4.1}.
Consider $a\not\in\Z$ and $A\in\CL^{a,k}(M,E)$. Furthermore,
choose $P\in\CL^{m,0}(M,E)$, $m>0$, self--adjoint elliptic whose
leading symbol is scalar and positive. 
Then by \myref{G1-3.10} the function $\Tr(AP^{-s})$ is regular at
$s=0$. Analogously to Definition \plref{S1-4.1} we put
\begin{eqnarray}
   \TR(A,P)&:=& \Tr(AP^{-s})|_{s=0}\\
    &=&\mbox{coefficient of}\, t^0\,\log^0 t\,\mbox{in the
asymptotic}\\
  &&\mbox{expansion of}\,\Tr(Ae^{-tP})\nonumber
\end{eqnarray}

\begin{theorem}{S5-5.1} $\TR(A):=\TR(A,P)$ is independent
of the particular $P$ chosen. $\TR$ defines a linear functional
on 
$$\bigcup_{a\in \C\setminus\Z,k\ge 0} \CL^{a,k}(M,E).$$
Furthermore,
\begin{enumerate}\renewcommand{\labelenumi}{{\rm (\roman{enumi})}}
\item $\TR\restr\CL^{a,k}(M,E)=\Tr_{\rmL^2}\restr\CL^{a,k}(M,E)$ if $a<-\dim M$.
\item $\TR([A,B])=0$ if $A\in\CL^{a,k}(M,E), B\in CL^{b,l}(M,E)$, $a+b\not\in\Z$.
\end{enumerate}
\end{theorem}

Thus, although the $\rmL^2$--trace cannot be extended as a trace on
$\CL^{*,*}$ there is an extension of the $\rmL^2$--trace to
non--integer order operators.

\proof 
It suffices to prove
the analogue of Proposition \plref{S1-4.3}. Obviously,
$\TR(A,P^\ga)=\TR(A,P)$. Thus w.l.o.g. we may assume $P$ to be
of integer order. The proof of loc. cit. shows
\begin{eqnarray*}
    \frac{d}{du}\Tr(e^{-tP_u})&\sim_{t\to 0+}&
      \sum_{j=0}^{\infty} \frac{(-1)^{j+1}}{(j+1)!}t^{j+1}\Tr(A (\nabla_P^j 
     \frac{d}{du} P_u)e^{-tP_u})\\
\Tr([A,B]e^{-tP})&\sim_{t\to 0+}&\sum_{j=1}^\infty
   \frac{(-1)^j}{j!} t^j 
   \Tr(A (\nabla_P^j B) e^{-tP}),
\end{eqnarray*}
and from \myref{G1-3.9} we immediately conclude that the coefficient
of $t^0$ in these expansions is $0$.\endproof

The {\sc Kontsevich--Vishik} trace has another interesting property
with respect to holomorphic families of operators. Furthermore,
the {\sc Kontsevich--Vishik} trace is given as the integral over
a canonical density. In order to generalize these facts
to our algebra we first introduce a regularized integral on the
space of symbols.

Consider $f\in\CS^{m,k}(\R^n)$. We write
\begin{equation}
    f=\sum_{j=0}^N \psi f_{m-j}+g,
\end{equation}
with $f_{m-j}\in\cp^{m-j,k}(\R^n), g\in
\CS^{m-N-1,k}(\R^n)$. $\psi\in\cinf{\R^n}$,
$\psi(\xi)=0$ if $|\xi|\le 1/4$ and $\psi(\xi)=1$ if $|\xi|\ge 1/2$.
In the sequel $\psi$ will be fixed. Then
\begin{equation}
    g(\xi)=O(|\xi|^{m-N}),\quad |\xi|\to\infty.
\end{equation}
This implies the asymptotic expansion
\begin{equation}
    \int_{|\xi|\le R}f(\xi)d\xi \sim_{R\to \infty}
    \sum_{j=0, j\not=m+n}^\infty p_{m+n-j}(\log R)\, R^{m-j+n}+
     p_0(\log R)\, R^0, 
    \mylabel{G5-5.7}
\end{equation}
with polynomials $p_\ga$ of degree
\begin{equation}
       \deg p_\ga\le\casetwo{k}{\ga\not=0,}{k+1}{\ga=0.}
\end{equation}
To see this we note that
\begin{equation}
    \int_{|\xi|\le R} g(\xi)d\xi=\int_{\R^n}g(\xi)d\xi
    +O(|\xi|^{m-N+n}),\quad |\xi|\to\infty,
\end{equation}
and splitting the integral over $\psi f_{m-j}$ into
$\int_{|\xi|\le 1}+\int_{1\le |\xi|\le R} $ we obtain
\begin{equation}
    \int_{1\le |\xi|\le R} \psi(\xi)f_{m-j}(\xi) d\xi
   =\sum_{l=0}^k \int_1^R\int_{|\xi|=1}f_{m-j,l}(\xi)d\xi\,
     r^{m-j+n-1}\, \log^l r dr,
   \mylabel{G5-5.10}
\end{equation}
which implies \myref{G5-5.7}.

We then define $\reginttext_{\R^n} f(\xi) d\xi$ to be the constant term in
the asymptotic expansion \myref{G5-5.7}, i.e.
\begin{equation}
\regint_{\R^n} f(\xi)d\xi:= \LIM_{R\to \infty} \int_{|\xi|\le R} f(\xi) d\xi
   :=p_0(0).
   \label{G5-5.11}
\end{equation}
Note that \myref{G5-5.10} implies that the coefficient of $R^0\,\log^{k+1}R$
in the expansion \myref{G5-5.7} equals
\begin{equation}
     \frac{\res_k(f_{-n})}{k+1}.
     \label{G5-5.12}
\end{equation}

\begin{prop}{S5-5.2} {\rm (cf. \cite[Lemma 2.1.4]{Les:OFTCSAM})}
Let $A\in {\rm GL}(n,\R)$ be a regular matrix
and $f\sim \sum\limits_{j\ge 0} f_{m-j} \in\CS^{m,k}(\R^n)$.
We have the transformation rule
$$\regint_{\R^n} f(A\xi) d\xi=
   |\det A|^{-1}\left( \regint_{\R^n} f(\xi)d\xi+
      \sum_{l=0}^{k}\frac{(-1)^{l+1}}{l+1}\int_{S^{n-1}} f_{-n,l}(\xi)
          \log^{l+1} |A^{-1}\xi| d\xi\right).$$
\end{prop}
\begin{Proof} It suffices to prove the proposition for $f\in \cinf{\R^n}$
with
$$f(\xi)=f(\xi/|\xi|)|\xi|^\ga \log^l|\xi|, \quad |\xi|\ge 1.$$
Then
\begin{equation}
\int_{|\xi|\le R} f(\xi) d\xi = \int_{|\xi|\le 1} f(\xi) d\xi
     + \int_1^R \int_{S^{n-1}} f(\xi) d\xi \,r^{n+\ga-1} \log^l r dr,
\end{equation}
and hence
\begin{equation}
   \regint_{\R^n} f(\xi) d\xi = \int_{|\xi|\le 1} f(\xi) d\xi+
     \displaycasetwo{\frac{(-1)^{l+1}l!}{(\ga+n)^{l+1}}
     \int_{S^{n-1}} f(\xi)d\xi}{\ga\not=-n,}{0}{\ga=-n.}
\end{equation}
On the other hand for $R$ large
\begin{eqnarray*}
  \int_{|\xi|\le R} f(A\xi) d\xi &=&  |\det A|^{-1}
    \int_{|A^{-1}\xi|\le R} f(\xi) d\xi\nonumber \\
    &=&  |\det A|^{-1} \left( \int_{|\xi|\le 1} f(\xi) d\xi +
      \int_{|\xi|\ge 1, |A^{-1}\xi|\le R}f(\xi) d\xi\right),
\end{eqnarray*}
\begin{eqnarray*}
    &&\int_{|\xi|\ge 1, |A^{-1}\xi|\le R} f(\xi)d\xi\\
   &=&\int_{S^{n-1}} f(\xi) \int_1^{R/|A^{-1}\xi|}
             r^{\ga+n-1} \log^l r dr d\xi\\
   &=&\left\{\renewcommand{\arraystretch}{2.5}\begin{array}{ll}\DST
  \int_{S^{n-1}} f(\xi) \Big(
   \big(\frac{R}{|A^{-1}\xi|}\big)^{\ga+n}
   \sum_{j=0}^l \frac{(-1)^{l-j}l!}{j!} \log^j (R/|A^{-1}\xi|)\; 
   (\ga+n)^{j-l-1}& \\
  \DST+\quad\quad
   \frac{(-1)^{l+1} l!}{(\ga+n)^{l+1}}\Big)d\xi,&\DST\alpha\not=-n,\\
  \DST \int_{S^{n-1}} f(\xi)
   \frac{1}{l+1} \log^{l+1} (R/|A^{-1}\xi|)d\xi ,  & \DST\alpha=-n,
      \end{array}\right.
\end{eqnarray*}
thus
\begin{equation}
   \LIM_{R\to\infty} \int_{|\xi|\ge 1,|A^{-1}\xi|\le R} f(\xi)d\xi=
   \left\{\renewcommand{\arraystretch}{2.5}\begin{array}{ll}
   \DST \frac{(-1)^{l+1} l!}{(\ga+n)^{l+1}}\int_{S^{n-1}}
   f(\xi)d\xi, & \ga\not=-n,\\
   \DST \frac{(-1)^{l+1}}{l+1}\int_{S^{n-1}}f(\xi)
   \log^{l+1} |A^{-1}\xi| d\xi,& \ga=-n,
   \end{array}\right.
\end{equation}
and we reach the conclusion. \end{Proof}

Next, we consider an open subset $U\subset \R^n$ and a pseudodifferential
operator $A\in\CL^{a,k}(U,\C^p,\C^q)$, $a\not\in\Z$, with amplitude
$a\in\CS^{a,k}(U,\R^n)\otimes{\rm Hom}(\C^p,\C^q)$, i.e.
\begin{equation}
   Au(x)=(2\pi)^{-n} \int_{\R^n} \int_U a(x,y,\xi) u(y) e^{i<x-y,\xi>}
      dy d\xi.
\end{equation}
For fixed $x$ we have $a(x,\cdot)\in\CS^{a,k}(\R^n)\otimes{\rm Hom}(\C^p,\C^q)$
and thus we may put
\newcommand{\KVden}{\omega_{\rm KV}}
\begin{equation}
    \KVden(A):=(2\pi)^{-n}\regint_{\R^n} a(x,x,\xi)d\xi\;|dx|.
\end{equation}

$\KVden$ has similar properties as the residue density:

\begin{lemma}{S5-5.3} We have
\begin{enumerate}\renewcommand{\labelenumi}{{\rm \arabic{enumi}.}}
\item $\KVden(A)\in\cinf{U,{\rm Hom}(\C^p,\C^q)\otimes|\Omega|}$.
\item $A\mapsto \KVden(A)$ is linear.
\item If $\varphi\in\cinfz{U}$ then $\KVden(\varphi A)=\varphi\KVden(A).$
\item If $\kappa:U\to V$ is a diffeomorphism then
  $$\kappa^*(\KVden(\kappa_*A))=\KVden(A).$$
\end{enumerate}
\end{lemma}
\proof 1.--3. are straightforward. To prove 4. we denote variables
in $U$ by $x,y$ and variables in $V$ by $\tilde x,\tilde y$.
$\kappa_*A$ has the amplitude function
\begin{equation}
     (\tilde x,\tilde y,\xi)\mapsto a(\kappa^{-1}\tilde x,\kappa^{-1}\tilde y,
   \phi(\tilde x,\tilde y)^{-1}\xi) \frac{|\det D\kappa^{-1}(\tilde x,\tilde
   y)|}{|\det \phi(\tilde x,\tilde y)|}
\end{equation}
(cf. \cite[Sec. 4.1, 4.2]{Shu:POST}), where $\phi(\tilde x,\tilde y)$ is
smooth with $\phi(\tilde x,\tilde x)=D\kappa^{-1}(\tilde x)^t$. Since
$a\not\in\Z$ the preceding proposition gives
\begin{equation}\begin{array}{rcl}
     \DST \KVden(\kappa_*A)&=&\DST (2\pi)^{-n}
               \regint_{\R^n} a(\kappa^{-1}\tilde x,\kappa^{-1}\tilde x,
   \phi(\tilde x,\tilde x)^{-1}\xi)d\xi\,|d\tilde x|\\[1.5em]
   &=&\DST(2\pi)^{-n} \regint_{\R^n} a(\kappa^{-1}\tilde x,\kappa^{-1}\tilde x,
   \xi)d\xi\,|\det D\kappa^{-1}(\tilde x)||d\tilde x|,
		\end{array}
\end{equation}
and thus
$$\kappa^*\KVden(\kappa_*A)=(2\pi)^{-n} \regint_{\R^n} a(x,x,\xi)d\xi\,
   |dx|=\KVden(A).\epformel$$

By 4. of the preceding lemma $\KVden$ can be defined globally on a
manifold. So, for $A\in\CL^{a,k}(M,E,F)$, $a\not\in\Z$, we obtain
a well--defined density
\begin{equation}
    \KVden(A)\in\cinf{M,{\rm Hom}(E,F)\otimes |\Omega|}.
   \label{G5-5.13}
\end{equation}
Furthermore, it is obvious that if $a<-\dim M$ and
$E=F$ then we have
\begin{equation}
   \Tr(A)=\int_M \tr_{E_x}(\KVden(A)(x)).
   \label{G5-5.14}
\end{equation}

Following {\sc Kontsevich--Vishik} \cite{KonVis:DEPO,KonVis:GDEO}
we now introduce holomorphic families.
Let $G\subset \C$ be a domain. A family of symbols $a(z)\in\CS^{z,k}(U,\R^n),
z\in G$, is called holomorphic if
\begin{equation}
   a(z)\sim \sum_{j=0}^\infty \psi a_{z-j},
\end{equation}
where
\begin{equation}
   a_{z-j}(z,x,\xi)=\sum_{l=0}^k a_{z-j,l}(z,x,\xi)\,\log^l|\xi|,
\end{equation}
$a_{z-j,l}(z,x,\cdot)\in\cp^{z,0}(\R^n)$, and
\begin{equation}
    z\mapsto a_{z-j}(z,\cdot,\cdot)
\end{equation}
is a holomorphic map into $\cinf{U\times\R^n}$. Furthermore,
\begin{equation}
     z\mapsto a(z)-\sum_{j=0}^N \psi a_{z-j}(z,\cdot,\cdot)
\end{equation}
is holomorphic with values in $C^{K(N)}(U\times\R^n), K(N)\to\infty$
as $N\to \infty$.

A family $A(z)\in\CL^{z,k}(M,E)$ is called holomorphic if in 
every chart $A(z)$ is given by a holomorphic amplitude $a(z)
\in\CS^{z,k}(U\times U,\R^n)\otimes {\rm End}(\C^p)$.

\begin{lemma}{S5-5.4} Let $f(z)\in\CS^{z,k}(\R^n)$, $z\in G$, be a holomorphic
family. Then the function
$$I(z):=\regint_{\R^n} f(z) dz$$
is meromorphic in $G$ with poles in $\Z\cap G$ of order $\le k+1$.
Furthermore,
$$ \Res_{k+1} I(z)|_{z=\nu}=(-1)^{k+1} k!\,\res_k(f_{-n}(\nu,\cdot)).$$
\end{lemma}
\proof Writing
\begin{equation}
   f(z)=\sum_{j=0}^N f_{z-j}+g(z)
\end{equation}
with $g(z,\xi)=O(|\xi|^{z-N}), |\xi|\to\infty$, we find that
\begin{equation}
    z\mapsto \regint_{\R^n} g(z,\xi)d\xi
\end{equation}
is holomorphic for $z\in G, \Re z<-n+N$.

Next consider a function of the form
\begin{equation}
      f(z,\xi)=\psi(\xi)f_{z-j,l}(z,\xi) \log^l|\xi|,\quad l\le k.
\end{equation}
Then $\int_{|\xi|\le 1} f(z,\xi)d\xi$ is holomorphic for $z\in G$ and
from \myref{G5-5.7} one derives
\begin{equation}
    \LIM_{R\to\infty}\int_{1\le |\xi|\le R} f_{z-j,l}(z,\xi)\log^l|\xi|d\xi
   =\frac{(-1)^{l+1}l!}{(z+n-j)^{l+1}}\int_{|\xi|=1} f_{z-j,l}(z,\xi)d\xi
\end{equation}
and we reach the conclusion.\endproof

\begin{prop}{S5-5.5} Let $A(z)\in \CL^{z,k}(M,E),z\in G$, be a holomorphic
family of operators. Then the function
    $$I(z):=\int_M\tr_{E_x}(\KVden(A(z))(x))$$
is meromorphic in $G$ with poles in $\Z\cap G$ of order $\le k+1$.
Moreover,
$$\Res_{k+1} I(z)|_{z=\nu}=\frac{(-1)^{k+1}}{k+1}\Res_{k}(A(\nu)).$$
\end{prop}
Note that $\Res_{k+1}$ on the left hand side means the coefficient of
$(z-\nu)^{-k-1}$ in the Laurent expansion (cf. \myref{G1-3.11})
while $\Res_{k}$ on the right
hand side means the noncommutative 
residue as defined in Definition \plref{S1-4.1}.
\proof This is a straightforward consequence of the preceding lemma
and Corollary \plref{S1-4.8}.\endproof

Now we are ready to state and prove the main result of this section,
which is the natural generalization of \cite[Sec. 3]{KonVis:DEPO}.

\begin{theorem}{S5-5.6} For $a\in\C\setminus\Z$ there exists a linear
functional
$$\TR:\CL^{a,k}(M,E)\longrightarrow \C$$
with the following properties:
\begin{enumerate}\renewcommand{\labelenumi}{{\rm (\roman{enumi})}}
\item For $A\in \CL^{a,k}(M,E)$ and any self--adjoint elliptic
$P\in\CL^{m,0}(M,E), m>0,$ with scalar and positive leading symbol
we have
$$\TR(A)=\Tr(AP^{-s})|_{s=0}.$$
\item $\TR(A)=\DST\int_M\tr_{E_x}(\KVden(A)(x)).$
\item $\TR\restr\CL^{a,k}(M,E)=\Tr_{\rmL^2}\restr\CL^{a,k}(M,E)$ if $a<-\dim M$.
\item $\TR([A,B])=0$ if $A\in \CL^{a,k}, B\in \CL^{b,l}, a+b\not\in\Z$.
\item If $A(z)\in\CL^{z,k}(M,E), z\in G$, is a holomorphic family
then $\TR(A(z))$ is meromorphic with poles in $z\in G\cap\{-\dim M+j\,|\,
j\in\Z_+\}$ of order $\le k+1$.
One has
$$\Res_{k+1} \TR(A(z))|_{z=\nu}=\frac{(-1)^{k+1}}{k+1}\Res_{k}(A(\nu)).$$
\end{enumerate}
\end{theorem}
\proof We take (i) as definition for $\TR$. Then it
only remains to prove (ii). Let $A\in\CL^{a,k}(M,E)$. Choose
$P\in \CL^{1,0}(M,E)$ self--adjoint elliptic with scalar and positive
leading symbol. We consider the family
\begin{equation}
     A(z):=AP^{-a+z}.
\end{equation}
Then $A(z)$, $z\in\C$, is a holomorphic family of operators. If
$z<-\dim M$ then $A(z)$ is trace class and hence by
\myref{G5-5.14} and Theorem \plref{S5-5.1}
\begin{equation}\begin{array}{rcl}
   \TR(A(z))&=& \Tr(A(z))=\Tr(AP^{z-a})\\
      &=&\DST \int_M \tr_{E_x}(\KVden(A(z))(x)).
		\end{array}
\end{equation}
Since left and right hand side extend meromorphically to $\C$ we find in
particular
$$\TR(A)=\TR(A(a))=\int_M\tr_{E_x}(\KVden(A)(x)).\epformel$$